# Structural Vibration Control of Offshore Wind Turbines Using Tuned Mass Damper Inerter in OpenFAST: Implementation, Validation, and Illustration

Hisham Tariq[1,3], Yuan Li[2], Flavia De Luca[3], and Agathoklis Giaralis[1]

[1] Department of Civil and Environmental Engineering, Khalifa University, Abu Dhabi, United Arab Emirates
[2] School of Electrical, Electronic and Mechanical Engineering, University of Bristol, Bristol, United Kingdom
[3] School of Civil, Aerospace and Design Engineering, University of Bristol, Bristol, United Kingdom

## Abstract

The tuned mass damper inerter (TMDI) is a passive vibration absorber widely considered in the literature to suppress structural vibrations by leveraging the inertance property of inerter devices. Previous studies have demonstrated TMDI effectiveness for vibration mitigation in wind turbines; however, these have relied on simplified structural models. Nevertheless, no study to date has implemented multi-directional TMDIs within aero-hydro-servo-elastic simulation tools such as OpenFAST to enable independent vibration control across multiple wind turbine components, including the tower, blades, and substructure. This paper addresses the above gap by implementing the TMDI in the OpenFAST. The governing equations of motion for TMDIs are derived and implemented within the Structural Control (StC) Module of OpenFAST, supporting arbitrary TMDI connectivity within each wind turbine component with independent multi-directional configurations. Verification against independent numerical models confirms the accuracy of the implementation across all supported configurations. An application study examines the IEA-15 MW reference wind turbine equipped with a tower-top TMDI under combined wind-wave loading conditions representative of normal power production and extreme idling scenarios. Optimal TMDI tuning is performed using a simplified two-degree-of-freedom model, whose modal properties are extracted from an OpenSeesPy finite element model of IEA-15 MW. Numerical results demonstrate that TMDIs with secondary mass up to 100 times smaller than a conventional TMD achieve matching or superior vibration suppression in the fore-aft and side-side directions, as quantified in terms of peak and standard deviation of both tower-top displacement and acceleration. Notably, the TMDI stroke is markedly reduced compared to the TMD, making it more suitable for integration within the spatially constrained wind turbine environment. Collectively, the foregoing results support the practical merit of lightweight TMDIs for wind turbine vibration control, and the developed OpenFAST implementation paves the way towards their performance-driven design under realistic loading conditions.

## 1. Introduction

Wind farms have emerged as a critical infrastructure in the global transition towards renewable energy, with offshore wind sector projecting to add 156 GW of new renewable energy generation capacity between 2025 and 2030 [1]. This is primarily achieved by favouring offshore wind turbines (OWTs) with larger rotor diameters and taller support structures offering access to higher and more consistent wind resources compared to onshore installations. In this regard, OWTs are becoming increasingly slender and flexible, rendering them particularly susceptible to excessive vibrations under combined environmental and operational excitations [2]. Indeed, wind turbines (WTs) are subjected to complex multi-directional vibrations encompassing tower deflections in fore-aft (FA) and side-side (SS) directions and blade out-of-plane (flapwise) and in-plane (edgewise) oscillations [3,4]. Additionally, floating support structures of OWTs experience platform rigid-body motions in surge, sway, heave, roll, pitch, and yaw [5,6]. These coupled dynamic responses may lead to accelerated fatigue damage, reduced energy generation efficiency, and ultimately, structural failure [2]. To this end, effective vibration suppression strategies are essential for next generation WTs [3,7], particularly as blade pitch control loses efficacy with increasing turbine power rating owing to excessive actuator force demands and blade-root stresses [8].

One of the most widely adopted strategies for WT vibration mitigation is the use of passive tuned mass dampers (TMDs) [9–14], comprising a secondary mass attached to the WT component (primary structure) via spring and damper elements tuned to a target vibration frequency. The TMD operates by transferring kinetic energy from the primary structure to the oscillating secondary mass and dissipating it through viscous damping. The practical merit of TMDs for WT applications has been demonstrated through their deployment in several commercial turbines [15,16]. However, the motion control effectiveness of the TMD is fundamentally constrained by two practical limitations: the requirement for considerable secondary mass (typically 1-5% of primary structure) [11,17] which imposes significant dead load, and prohibitively large oscillation strokes demanding substantial clearance space [16,17]. These stroke requirements become particularly critical for flexible WTs with low natural frequencies, posing significant challenges for integration within the spatially restricted WT system and increasing upfront costs.

To address these TMD limitations, the inerter, introduced in [18] as a two-terminal mechanical element resisting relative acceleration through its inertance property, has been leveraged to enhance passive TMDs in recent years giving rise to inerter enhanced vibration absorbers (IVAs). The inerter achieves inertial amplification effect by transforming translational motion into high-speed rotational motion, enabling inertance values several orders of magnitude larger than the physical device mass. In the context of WTs, single point IVAs have been extensively pursued in literature [19–28]. Here 'single-point' refers to IVA configurations that connect the primary structure at a single node such that the inerter element of inerter based network is placed in parallel with the conventional spring-damper assembly of the TMD. Fig. 1a illustrates the simplest such configuration and Qian and Liao [22] developed an analytical model for this IVA topology coupled with wind turbine towers and evaluated its vortex-induced vibration mitigation performance in the NREL 5MW onshore WT.

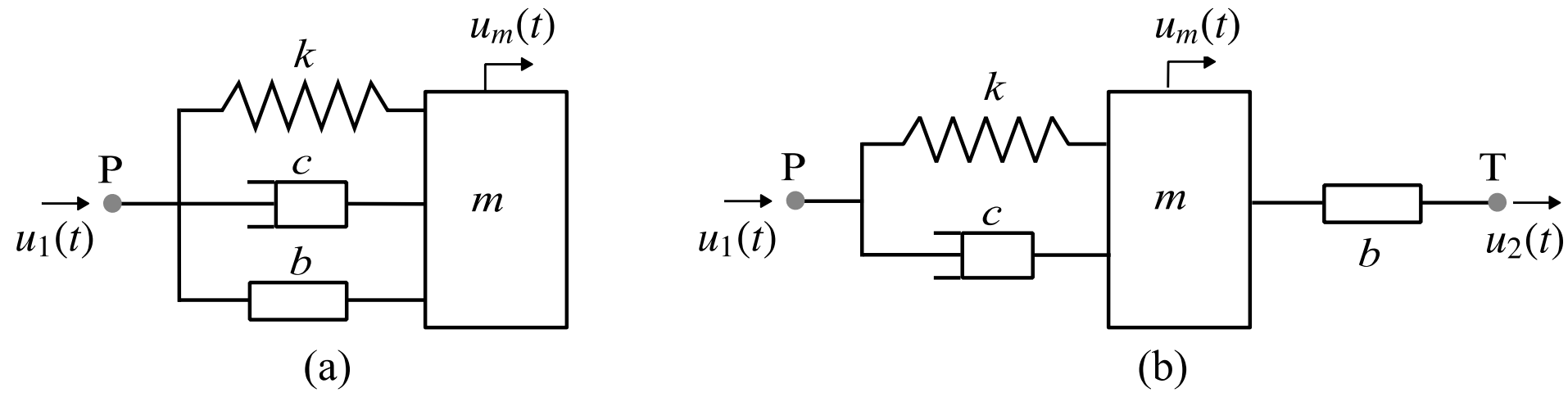

**Fig. 1** Schematic of two IVA topologies: (a) single-point IVA with spring, damper, and inerter connected in parallel between structural node P and the absorber mass; and (b) TMDI with spring and damper connecting the absorber mass to structural node P and the inerter connecting it to structural node T.

While the above single-point IVAs have demonstrated promising vibration-suppression performance, achieving substantial response reductions typically relies on relatively large inertance values. This is attributed to the fact that the inerter-based network acts solely on the relative motion between the absorber mass and a single structural attachment point. In this respect, two-point topologies that exploit relative motion between two structural locations have been pursued. Among these, the tuned mass damper inerter (TMDI), introduced by Marian and Giaralis [29], features a series connection of an inerter element with the TMD (see Fig. 1b). In this configuration, the TMDI is attached to the WT component at two locations: the secondary mass is connected at one structural location via spring and damper elements acting in parallel, while the inerter connects this secondary mass to a different structural location, thereby achieving enhanced motion control performance [30,31]. Specifically, Sarkar and Fitzgerald [31] demonstrated the effectiveness of TMDI for vibration control in FOWTs, wherein one terminal was attached to the nacelle and the other to the tower. This implementation achieved remarkable performance improvements under both normal operational and extreme wind-wave coupled loading conditions. Subsequently, Sarkar and Fitzgerald [32] investigated the tuned mass damper fluid-inerter (TMDFI) for vibration control of FOWT towers, demonstrating that the fluid-inerter achieves performance equivalent to the ideal mechanical inerter, with the TMDFI outperforming conventional TMDs considerably. More recently, Fitzgerald et al. [33] examined TMDI implementation for enhancing structural reliability of FOWT towers under misaligned wind-wave loading conditions. Additionally, McAuliffe et al. [34] demonstrated the superior fatigue life enhancement of FOWT towers achieved by TMDI compared to conventional TMD. Furthermore, Peng and Chang [35] employed reliability-based design optimization (RBDO) to demonstrate TMDI efficacy in vibration mitigation of FOWTs. Extending TMDI applications to blade vibration control, Zhang and Fitzgerald [36] investigated TMDI implementation for suppressing edgewise vibrations in blades, demonstrating considerable reductions in damper stroke compared to classical TMDs. Further investigations by [37,38] examined TMDI for seismic vibration control and energy harvesting applications. Additionally, Shi and Zhu [39] employed a geometrically nonlinear TMDI formulation to assess dynamic response reduction in wind turbines, while Alotta et al. [40] studied the tuned inerter damper (TID), a TMDI configuration with zero secondary mass, for seismic response mitigation of land-based WTs. Collectively, these studies demonstrate the versatility and effectiveness of TMDI configurations for vibration suppression across diverse wind turbine applications and loading scenarios.

Despite these advances, most TMDI applications for WTs reported to date have relied on either lumped-parameter reduced-order models [37–39] or multibody-dynamics formulations [31–34] benchmarked against FAST (a predecessor to OpenFAST) under idealized steady-state wind and still-water conditions. Of particular relevance, Zhang et al. [30] implemented a unidirectional tower top TMDI in FAST, configured to act solely in the FA direction. Nevertheless, the existing StC module of OpenFAST v5.0 does not support TMDIs. Rather, it only supports single attachment point configurations without inerter elements [41], thereby preventing the design and performance assessment of TMDIs across multiple degrees of freedom within a high-fidelity aero-hydro-servo-elastic simulation environment under realistic operational and environmental loading conditions. To the best of the authors' knowledge, multi-directional TMDIs with arbitrary intra-component TMDI connectivity across the tower, nacelle, substructure, and blades have not been implemented within the OpenFAST StC module to date.

To this end, this paper extends the OpenFAST StC module to support multi-directional TMDIs with arbitrary intra-component connectivity across the tower, nacelle, substructure, and blades. First, the equations of motion for the TMDI are derived using multibody-dynamics principles, consistent with the modular structural-dynamics framework of OpenFAST. Next, these equations are integrated into the StC module through source code modification, in accordance with OpenFAST's code architecture and programming conventions. The implementation is then verified for accuracy against independent numerical models across tower-, substructure-, and blade-based TMDI configurations.

The practical merit of the TMDI implementation in OpenFAST is further demonstrated through a case study on tower-top vibration mitigation using an optimally designed bidirectional TMDI under combined wind-wave loading. The case study employs the IEA-15 MW offshore reference wind turbine as a benchmark for contemporary large-scale OWT design. For this purpose, optimal TMDI parameters along the FA direction and SS direction are computed using a two-degree-of-freedom system whose modal properties are extracted from an OpenSeesPy finite element model of IEA-15 MW. Design load cases representing normal operating conditions at rated wind speed and parked conditions in case extreme storm are considered. Performance of TMD(I) equipped IEA 15MW is quantified in terms of tower top displacements, tower top accelerations, TMD(I) strokes and inerter forces in FA and SS direction.

The remainder of the paper is structured as follows: Section 2 describes the TMDI modelling approach in OpenFAST, including the derivation of the TMDI equations of motion and their implementation in the StC module. Section 3 presents the numerical validation of the TMDI implementation in OpenFAST. Section 4 carries out the assessment of the TMD(I)-equipped IEA 15MW wind turbine, and Section 5 summarizes the concluding remarks.

## 2. Proposed modelling approach of the tuned mass damper inerter in OpenFAST

In this section, the derivation of the equations of motion of the tuned mass damper inerter in OpenFAST three-dimensional vector-space notation is presented and the integration of these equations in state-space, following the multi-body dynamics formulation used in OpenFAST, is discussed. The presentation begins by a brief overview of the modular architecture of OpenFAST which facilitates the implementation of TMDI in different OWT components (e.g. tower, blades, substructure).

### 2.1 Modular Architecture of OpenFAST and Structural Control Module Extensions

Originating as an acronym from "Fatigue, Aerodynamics, Structures, and Turbulence", OpenFAST has emerged as the *de facto* open-source simulation platform for detailed aero-hydro-servo-elastic time-domain analyses, capturing faithfully the dynamic response of OWTs under realistic metocean (wind and wave) conditions. This is achieved by a modular multi-physics/domain architecture in which each different domain, management by a dedicated module, is interacting with all other domains through a common coupler referred to as the glue code. The glue code coordinates module initialization, manages data exchange across motion and load meshes, and advances the coupled simulation in time. The principal modules instantiating the OpenFAST modular architecture are InflowWind, AeroDyn, SeaState, HydroDyn, ElastoDyn, SubDyn, and ServoDyn. These modules collectively govern environmental excitation, aerodynamic and hydrodynamic loading, structural dynamics, substructure response, and control system behaviour of the OWT system, as illustrated in Fig. 2.

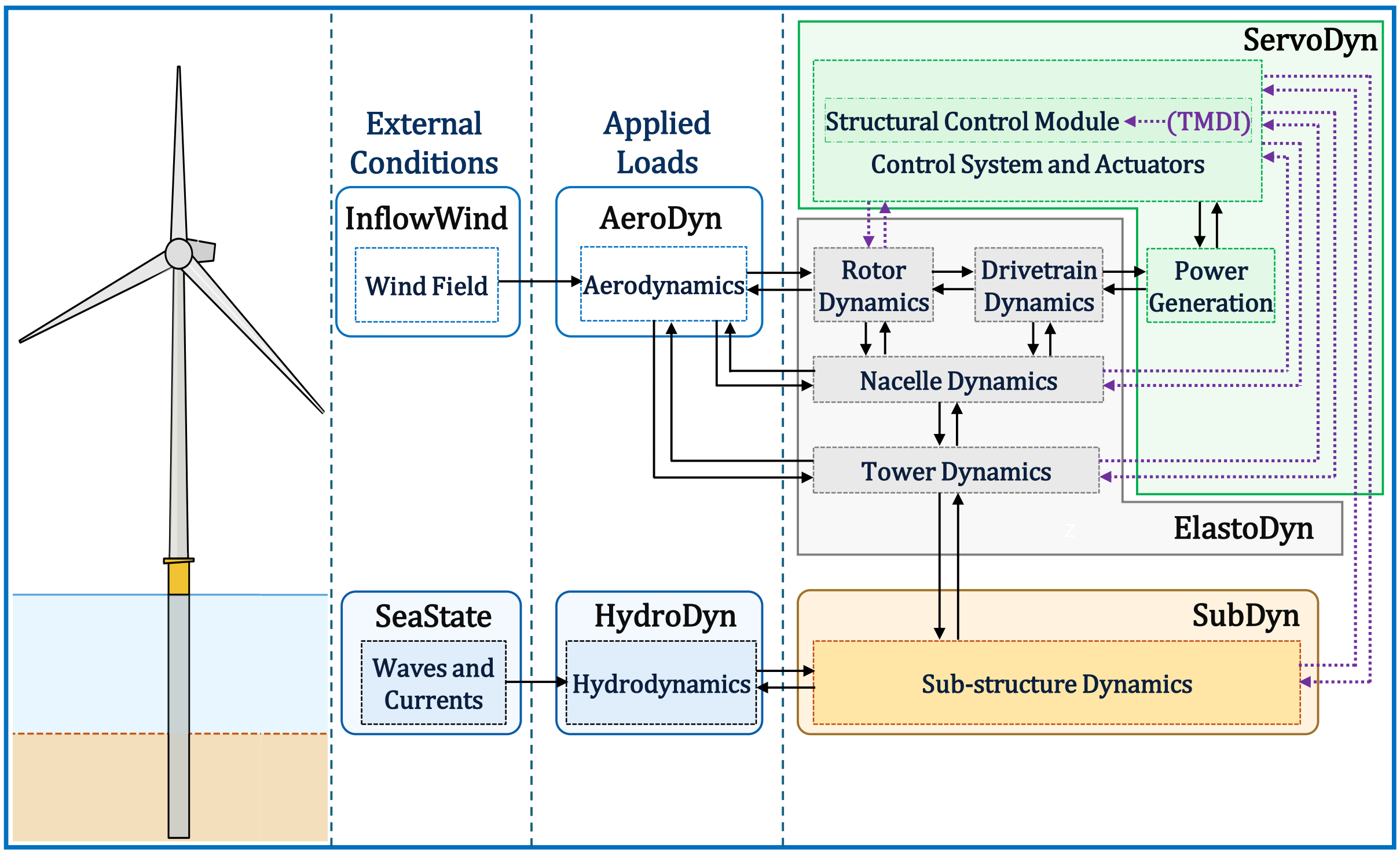

**Fig. 2** Modular architecture of OpenFAST for OWT simulation. Dotted arrows indicate data exchange pathways modified in the present work to support tuned mass damper inerter (TMDI) modelling.

Of most relevance to the present work is the structural control (StC) sub-module, which is part of the ServoDyn module as shown in Fig.2, used by the current and all the previous versions of OpenFAST to model various standard single-point dynamic vibration absorbers (DVAs), including TMDs, and account for DVA dynamic interaction with the OWT. Herein, the StC sub-module is extended to allow for inerter element modelling as well as for the additional kinematics and load coupling required by the two-point tuned mass damper inerter (TMDI) absorber in Fig.1(b), operating independently along three orthogonal directions: x, y, and z. The former extension is achieved by hard-coding into StC sub-module the TMDI equations of motion derived as detailed in subsequent subsections, while the latter extension is achieved by appropriately modifying the data exchange pathways of OpenFAST glue code indicated by dotted arrows in Fig. 2. Importantly, the modular architecture of OpenFAST enables the implementation of the TMDI, upon implementation of the herein developed extensions, in any OWT component, including the tower, the blades, or the substructure. For illustration, Fig. 3 presents possible TMDI implementations at different components, in which point P denotes the attachment location of the spring-damper assembly on the component and point T denotes the attachment of the inerter.

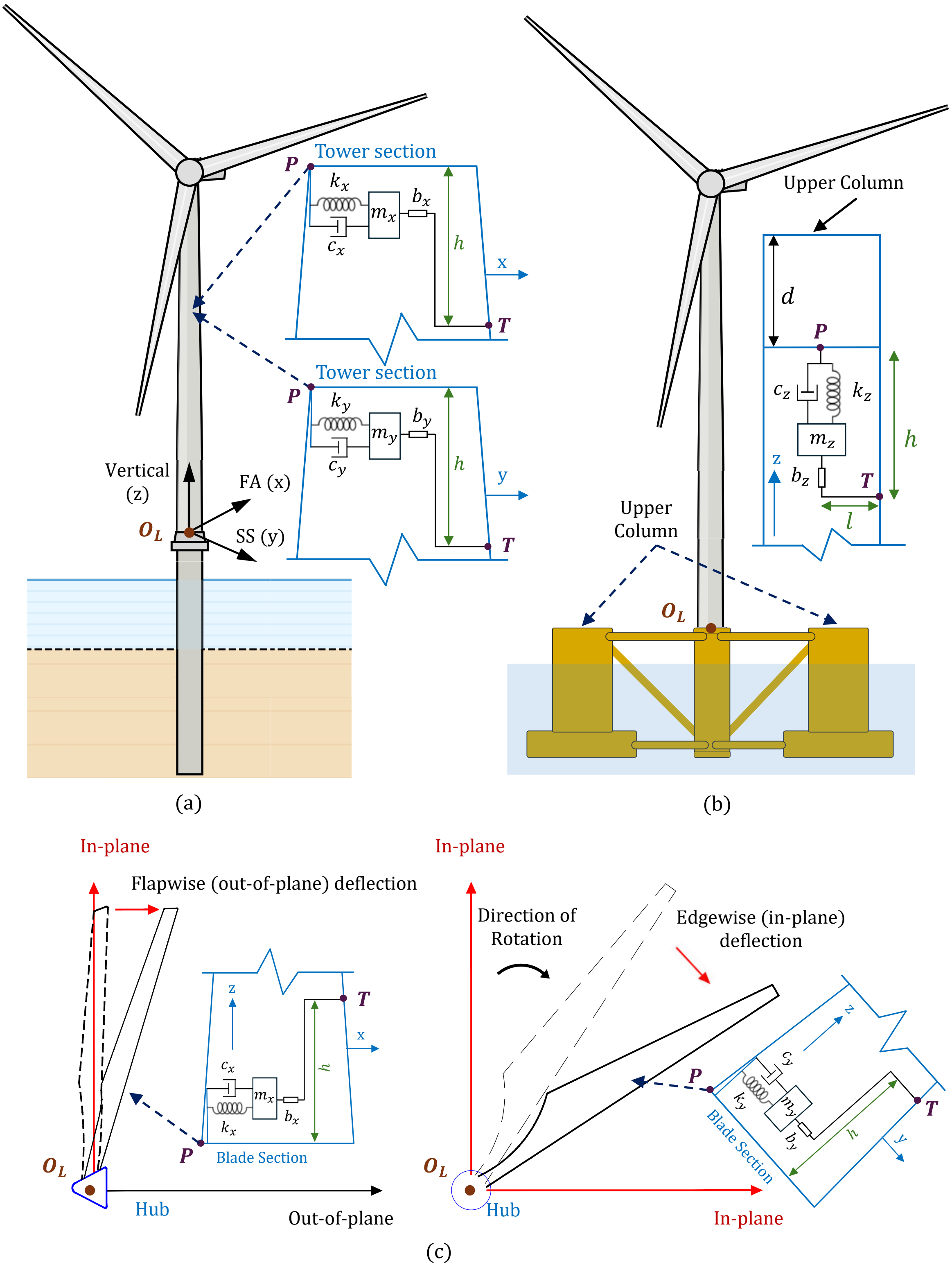

**Fig. 3** Schematic of TMDI configurations within various wind turbine components: (a) bi-directional TMDI at the tower top acting along the fore-aft (x) and side-side (y) directions; (b) three unidirectional TMDIs, one per upper column of the semi-submersible floating substructure, each acting along the vertical (z) direction; (c) bi-directional TMDIs at the blade roots, each acting along the flapwise (x) and edgewise (y) directions. The mass, inertance, stiffness and damping of TMDI acting along local direction $i$ ($i \in \{x, y, z\}$) are denoted by $m_i, b_i, k_i$, and $c_i$, respectively.

### 2.2 Derivation of multi body dynamics equations of TMDI mass

To derive the TMDI equations of motion and to specify the two-point kinematic and load coupling between the TMDI and its host component, two different coordinate orientations, native to OpenFAST, need to be followed. The global orientation, denoted by subscript $G$, has its origin at the undeflected tower centreline at ground level for onshore WTs and at mean sea level for offshore WTs. The local orientation, denoted by subscript $L$, is component-specific. For towers, the origin of local orientation is at the tower base, for blades is at the blade root, and for substructures is at the so-called platform reference point, defined differently for each type of substructrure (e.g. monopiles, jackets, floating platforms). The relationship between these two orientations is defined through the direction cosine matrix $R_{L/G}$, which is a $3 \times 3$ rotation matrix that transforms vectors from the global orientation $G$ into the local component orientation $L$. The inverse transformation from local to global is obtained by taking the transpose of $R_{L/G}$, denoted by $R^T_{L/G}$. The direction cosine matrix is computed by the structural module ElastoDyn or SubDyn, depending on the host component in which the TMDI is installed, and is subsequently passed through the glue code to the StC sub-module of ServoDyn as an input. Further to the above orientations, two reference frames are defined in OpenFAST as illustrated in Fig. 4: (i) an inertial (non-rotating) global reference frame with origin at point O, coinciding with the origin of the global orientation, and (ii) a non-inertial (rotating) local reference frame with origin at point P, associated with the component hosting the TMDI at rest. Note that point P may be distinct from the origin of the local orientation of the component, depending on the attachment location of the TMDI spring–damper assembly within the component.

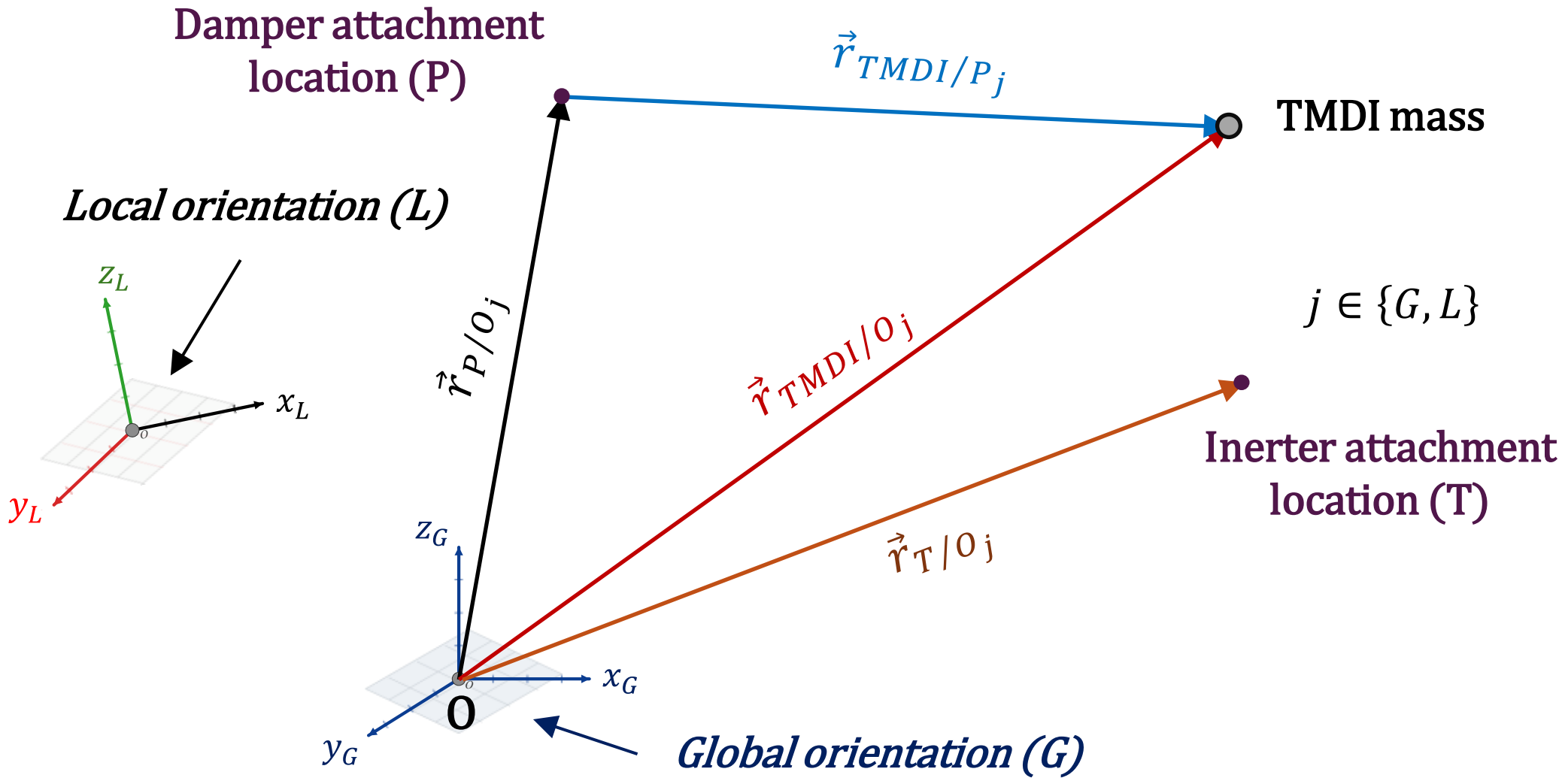

**Fig. 4** Schematic of the global (G) and local (L) orientations, the associated reference frame origins $O$ and $P$, and the position vectors to the TMDI mass and inerter attachment location used to define the TMDI connectivity within the StC module of OpenFAST.

The TMDI mass is modelled within OpenFAST as a particle possessing translational inertial properties (Fig. 4). To establish the governing equations of motion, the following assumptions are adopted. First, the inerter is idealized as a massless linear element, as defined by Smith [18]. Second, the inerter is assumed to act solely along a translational degree of freedom, with its rotational contributions to the TMDI dynamics neglected. Finally, the two terminals of the inerter element are assumed to be rigidly connected to the TMDI mass and to an arbitrary point T within the component.

Under these modelling assumptions the position vector of the TMDI mass relative to point P, expressed in the global orientation $\vec{r}_{TMDI/P_G}$ and the local component orientation $\vec{r}_{TMDI/P_L}$, can be written as

$$\begin{aligned} \vec{r}_{TMDI/P_G} &= \vec{r}_{TMDI/O_G} - \vec{r}_{P/O_G}, \\ \vec{r}_{TMDI/P_L} &= \vec{r}_{TMDI/O_L} - \vec{r}_{P/O_L}, \end{aligned} \tag{1}$$

respectively. In the above expressions, $\vec{r}_{TMDI/O_G}$ and $\vec{r}_{P/O_G}$ denote the position vectors of the TMDI mass and point P relative to O in the global orientation, respectively, while $\vec{r}_{TMDI/P_L}$ and $\vec{r}_{P/O_L}$ denote the position vectors of the TMDI mass and point P relative to O in the local orientation, respectively. The velocity vector of TMDI relative to P in the local orientation can be computed using the transport theorem as [42,43]

$$\dot{\vec{r}}_{TMDI/P_L} = \dot{\vec{r}}_{TMDI/O_L} - \dot{\vec{r}}_{P/O_L} - \vec{\omega}_{L/O_L} \times \vec{r}_{TMDI/P_L}, \tag{2}$$

where a dot over a symbol signifies differentiation with time and $\vec{\omega}_{L/O_L} = [\dot{\theta} \;\; \dot{\phi} \;\; \dot{\psi}]^T_{L/O_L}$ is the angular velocity of the component with respect to O in the local orientation. Further differentiation of Eq. 2 with respect to time, followed by application of the transport theorem, yields the acceleration of the TMDI mass relative to P in the local orientation as

$$\begin{aligned} \ddot{\vec{r}}_{TMDI/P_L} = &\frac{1}{m}\vec{F}_{TMDI/O_L} - \ddot{\vec{r}}_{P/O_L} - \vec{\alpha}_{L/O_L} \times \vec{r}_{TMDI/P_L} - \vec{\omega}_{L/O_L} \times (\vec{\omega}_{L/O_L} \times \vec{r}_{TMDI/P_L}) \\ &-2\vec{\omega}_{L/O_L} \times \dot{\vec{r}}_{TMDI/P_L}, \end{aligned} \tag{3}$$

where $\vec{F}_{TMDI/O_L}$ is the total force acting on TMDI relative to O, $m$ is the TMDI mass, and $\ddot{\vec{r}}_{P/O_L}$ is the translational acceleration of point P with respect to O. Further, in Eq.(3), the last three terms on the right-hand side account for rotational effects. Specifically, first, the term $\vec{\alpha}_{L/O_L} \times \vec{r}_{TMDI/P_L}$ is the tangential acceleration, which arises when the point P undergoes angular acceleration $\vec{\alpha}_{L/O_L}$. Next, the term $\vec{\omega}_{L/O_L} \times (\vec{\omega}_{L/O_L} \times \vec{r}_{TMDI/P_L})$ is the centripetal acceleration directed towards the instantaneous axis of rotation, whose direction is $\vec{\omega}_{L/O_L}/|\vec{\omega}_{L/O_L}|$, with magnitude proportional to the square of the angular velocity. Lastly, the term $2\vec{\omega}_{L/O_L} \times$

$\dot{\vec{r}}_{TMDI/P_L}$ accounts for Coriolis acceleration, which emerges when the TMDI mass moves within the rotating frame and acts perpendicular to both the angular velocity vector $\vec{\omega}_{L/O_L}$ and the relative velocity $\dot{\vec{r}}_{TMDI/P_L}$. At this junction, Eq. (3) needs to be particularized for each different possible direction of action of the TMDI, along the local x, y, and z directions, as detailed in the following three sub-sections. This is necessitated by the fact that the relative position and velocity vectors depend on the TMDI direction of action, which further makes the associated tangential, centripetal, and Coriolis contributions to be direction-dependent.

### 2.3 TMDI equations of motion along x direction

First, the case of TMDI mass constrained to move only along the local x-axis, denoted by $TMDI_X$, is elaborated. Thus, only the $TMDI_X$ position, velocity, and acceleration components relative to P, namely, $x_{TMDI_X/P_L}$, $\dot{x}_{TMDI_X/P_L}$, and $\ddot{x}_{TMDI_X/P_L}$, respectively, are non-zero. By application of force equilibrium, the resultant force vector acting on TMDI relative to the origin, $\vec{F}_{TMDI_X/O_L}$, is written as

$$\vec{F}_{TMDI_X/O_L} = \begin{bmatrix} -b_x(\ddot{x}_{TMDI_X/P_L} + \ddot{x}_{P/O_L} - \ddot{x}_{T/O_L}) - c_x\dot{x}_{TMDI_X/P_L} - k_x x_{TMDI_X/P_L} + m_x a_{G_X/O_L} + F_{ext_X} + F_{StopFrc_X} \\ F_{Y_{TMDI_X/O_L}} + m_x a_{G_Y/O_L} \\ F_{Z_{TMDI_X/O_L}} + m_x a_{G_Z/O_L} \end{bmatrix}. \quad (4)$$

In Eq.(4), $a_{G_X/O_L}$, $a_{G_Y/O_L}$ and $a_{G_Z/O_L}$ denote the components of gravitational acceleration along x, y, and z directions, respectively, and $\ddot{x}_{T/O_L}$ denotes the acceleration of point T along the x direction in the local orientation. Further, the term $F_{ext_X}$ represents the external control force acting on the $TMDI_X$ mass along the x-direction, which enables active vibration suppression when the system operates in active mode, while in passive mode this force is zero. The stop spring force, $F_{StopFrc_X}$, activates when the TMDI mass displacement relative to P exceeds a prescribed limit. The constraint reactions forces acting on $TMDI_X$ mass. $F_{Y_{TMDI_X/O_L}}$ and $F_{Z_{TMDI_X/O_L}}$ prevent the motion of $TMDI_X$ mass along Y and Z directions caused from the centripetal, tangential, and Coriolis accelerations, thereby enforcing the unidirectional constraint. For the $TMDI_X$ case, Eq. (3) reduces to

$$\begin{bmatrix} \ddot{x}_{TMDI_X/P_L} \\ 0 \\ 0 \end{bmatrix} = \frac{1}{m_x}\vec{F}_{TMDI/O_L} - \ddot{\vec{r}}_{P/O_L} - \vec{\alpha}_{L/O_L} \times \vec{r}_{TMDI/P_L} - \vec{\omega}_{L/O_L} \times (\vec{\omega}_{L/O_L} \times \vec{r}_{TMDI/P_L}) - 2\vec{\omega}_{L/O_L} \times \dot{\vec{r}}_{TMDI/P_L}, \quad (5)$$

in which the rotational acceleration terms can be expressed as

$$\begin{aligned} &\vec{\alpha}_{L/O_L} \times \vec{r}_{TMDI/P_L} = x_{TMDI_X/P_L}\begin{bmatrix}0 & \ddot{\psi}_{L/O_L} & \ddot{\phi}_{L/O_L}\end{bmatrix}^T, \\ &\vec{\omega}_{L/O_L} \times (\vec{\omega}_{L/O_L} \times \vec{r}_{TMDI/P_L}) = x_{TMDI_X/P_L}\begin{bmatrix}-(\dot{\phi}^2_{L/O_L}+\dot{\psi}^2_{L/O_L}) & \dot{\theta}_{L/O_L}\dot{\phi}_{L/O_L} & \dot{\theta}_{L/O_L}\dot{\psi}_{L/O_L}\end{bmatrix}^T, \end{aligned} \quad (6)$$

$$2\vec{\omega}_{L/O_L} \times \dot{\vec{r}}_{TMDI/P_L} = \dot{x}_{TMDI_X/P_L}\left[0 \quad 2\dot{\psi}_{L/O_L} \quad -2\dot{\phi}_{L/O_L}\right]^T.$$

Substituting Eqs. (4) and (6) into Eq. (5) and projecting on the local x direction yields the governing equation of motion for the TMDI$_X$ case as

$$\begin{aligned}\ddot{x}_{TMDI_X/P_L} = &\left[\frac{m_x\left((\dot{\phi}^2_{L/O_L}+\dot{\psi}^2_{L/O_L})-k_x\right)}{m_x+b_x}\right]x_{TMDI_X/P_L} - \frac{c_x}{m_x+b_x}\dot{x}_{TMDI_X/P_L} - \frac{m_x}{m_x+b_x}\ddot{x}_{P/O_L} \\ &-\frac{b_x}{m_x+b_x}(\ddot{x}_{P/O_L}-\ddot{x}_{T/O_L}) + \frac{m_x}{m_x+b_x}a_{G_X/O_L} + \frac{1}{m_x+b_x}(F_{ext_x}+F_{StopFrc_X}),\end{aligned} \tag{7}$$

while the constraint reaction forces follow from the $y$ and $z$ components as

$$\begin{aligned}F_{Y_{TMDI_X/O_L}} &= m_x\left[-a_{G_Y/O_L} + \ddot{y}_{P/O_L} + (\ddot{\psi}_{L/O_L}+\dot{\theta}_{L/O_L}\dot{\phi}_{L/O_L})x_{TMDI_X/P_L} + 2\dot{\psi}_{L/O_L}\dot{x}_{TMDI_X/P_L}\right], \\ F_{Z_{TMDI_X/O_L}} &= m_x\left[-a_{G_Z/O_L} + \ddot{z}_{P/O_L} - (\ddot{\phi}_{L/O_L}-\dot{\theta}_{L/O_L}\dot{\psi}_{L/O_L})x_{TMDI_X/P_L} - 2\dot{\phi}_{L/O_L}\dot{x}_{TMDI_X/P_L}\right].\end{aligned} \tag{8}$$

### 2.4 TMDI equations of motion along y direction

In this sub-section the case of TMDI mass constrained to move only along the local y-axis, denoted by TMDI$_Y$, is considered. The non-zero TMDI position, velocity, and acceleration components relative to P are $y_{TMDI_Y/P_L}$, $\dot{y}_{TMDI_Y/P_L}$, and $\ddot{y}_{TMDI_Y/P_L}$respectively. By application of force equilibrium, the resultant force vector acting on TMDI relative to origin, $\vec{F}_{TMDI_Y/O_L}$, is given as

$$\vec{F}_{TMDI_Y/O_L} = \begin{bmatrix} F_{X_{TMDI_Y/O_L}} + m_y a_{G_X/O_L} \\ -b_y(\ddot{y}_{TMDI_Y/P_L} + \ddot{y}_{P/O_L} - \ddot{y}_{T/O_L}) - c_y\dot{y}_{TMDI_Y/P_L} - k_y y_{TMDI_Y/P_L} + m_y a_{G_Y/O_L} + F_{ext_y} + F_{StopFrc_Y} \\ F_{Z_{TMDI_Y/O_L}} + m_y a_{G_Z/O_L}\end{bmatrix}, \tag{9}$$

where $\ddot{y}_{T/O_L}$ denotes the acceleration of point T along the y direction in the local orientation, and the terms $F_{ext_Y}$, $F_{StopFrc_Y}$, $F_{X_{TMDI_Y/O_L}}$, and $F_{Z_{TMDI_Y/O_L}}$are analogous to the TMDI$_X$ case previously discussed. For the TMDI$_Y$ case, Eq. (3) reduces to

$$\begin{aligned}\begin{bmatrix}0 \\ \ddot{y}_{TMDI_Y/P_L} \\ 0\end{bmatrix} = &\frac{1}{m_y}\vec{F}_{TMDI/O_L} - \ddot{\vec{r}}_{P/O_L} - \vec{\alpha}_{L/O_L} \times \vec{r}_{TMDI/P_L} - \vec{\omega}_{L/O_L} \times (\vec{\omega}_{L/O_L} \times \vec{r}_{TMDI/P_L}) \\ &- 2\vec{\omega}_{L/O_L} \times \dot{\vec{r}}_{TMDI/P_L},\end{aligned} \tag{10}$$

in which the rotational acceleration terms are

$$\begin{aligned}\vec{\alpha}_{L/O_L} \times \vec{r}_{TMDI/P_L} &= y_{TMDI_Y/P_L}\left[-\ddot{\psi}_{L/O_L} \quad 0 \quad \ddot{\theta}_{L/O_L}\right]^T, \\ \vec{\omega}_{L/O_L} \times (\vec{\omega}_{L/O_L} \times \vec{r}_{TMDI/P_L}) &= y_{TMDI_Y/P_L}\left[\dot{\theta}_{L/O_L}\dot{\phi}_{L/O_L} \quad -(\dot{\theta}^2_{L/O_L}+\dot{\psi}^2_{L/O_L}) \quad \dot{\phi}_{L/O_L}\dot{\psi}_{L/O_L}\right]^T, \\ 2\vec{\omega}_{L/O_L} \times \dot{\vec{r}}_{TMDI/P_L} &= \dot{y}_{TMDI_Y/P_L}\left[-2\dot{\psi}_{L/O_L} \quad 0 \quad 2\dot{\theta}_{L/O_L}\right]^T.\end{aligned} \tag{11}$$

Substituting Eqs. (9) and (11) into Eq. (10) and projecting on the local y direction yields the governing equation of motion for TMDI$_Y$ case as

$$\ddot{y}_{TMDI_Y/P_L} = \left[\frac{m_y\left(\left(\dot{\theta}^2_{L/O_L}+\dot{\psi}^2_{L/O_L}\right)-k_y\right)}{m_y+b_y}\right] y_{TMDI_Y/P_L} - \frac{c_y}{m_y+b_y}\dot{y}_{TMDI_Y/P_L} - \frac{m_y}{m_y+b_y}\ddot{y}_{P/O_L}$$
$$-\frac{b_y}{m_y+b_y}\left(\ddot{y}_{P/O_L}-\ddot{y}_{T/O_L}\right) + \frac{m_y}{m_y+b_y}a_{G_Y/O_L} + \frac{1}{m_y+b_y}\left(F_{ext_Y}+F_{StopFrc_Y}\right), \quad (12)$$

with constraint reactions

$$F_{X_{TMDI_Y/O_L}} = m_y\left[-a_{G_X/O_L} + \ddot{x}_{P/O_L} - (\ddot{\psi}_{L/O_L}-\dot{\theta}_{L/O_L}\dot{\phi}_{L/O_L})y_{TMDI_Y/P_L} + 2\dot{\psi}_{L/O_L}\dot{y}_{TMDI_Y/P_L}\right],$$
$$F_{Z_{TMDI_Y/O_L}} = m_y\left[-a_{G_Z/O_L} + \ddot{z}_{P/O_L} + (\ddot{\theta}_{L/O_L}+\dot{\phi}_{L/O_L}\dot{\psi}_{L/O_L})y_{TMDI_Y/P_L} + 2\dot{\theta}_{L/O_L}\dot{y}_{TMDI_Y/P_L}\right]. \quad (13)$$

## 2.5 TMDI equations of motion along z direction

Lastly, the case of TMDI mass constrained to move only along the local z-axis, denoted by $TMDI_Z$, is addressed. The non-zero TMDI position, velocity, and acceleration components relative to P are $z_{TMDI_Z/P_L}$, $\dot{z}_{TMDI_Z/P_L}$, and $\ddot{z}_{TMDI_Z/P_L}$, respectively. By application of force equilibrium, the resultant force vector acting on TMDI relative to the origin, $\vec{F}_{TMDI_Z/O_L}$, is given as

$$\vec{F}_{TMDI_Z/O_L} = \begin{bmatrix} F_{X_{TMDI_Z/O_L}} + m_z a_{G_X/O_L} \\ F_{Y_{TMDI_Z/O_L}} + m_z a_{G_Y/O_L} \\ -b_z(\ddot{z}_{TMDI_Z/P_L} + \ddot{z}_{P/O_L} - \ddot{z}_{T/O_L}) - c_z\dot{z}_{TMDI_Z/P_L} - k_z z_{TMDI_Z/P_L} + m_z a_{G_Z/O_L} + F_{ext_Z} + F_{StopFrc_Z} + F_{Z_{PreLoad}} \end{bmatrix}. \quad (14)$$

where $\ddot{z}_{T/O_L}$ denotes the acceleration of point T along the z direction in the local orientation, $F_{Z_{PreLoad}}$ compensates for the static deflection of the spring connecting the component and $TMDI_Z$ mass due to gravity, and the terms $F_{ext_Z}$, $F_{StopFrc_Z}$, $F_{X_{TMDI_Z/O_L}}$, and $F_{Y_{TMDI_Z/O_L}}$ are analogous to the $TMDI_X$ and $TMDI_Y$ cases previously discussed. For the $TMDI_Z$ case, Eq. (3) reduces to

$$\begin{bmatrix} 0 \\ 0 \\ \ddot{z}_{TMDI_Z/P_L} \end{bmatrix} = \frac{1}{m_z}\vec{F}_{TMDI/O_L} - \ddot{\vec{r}}_{P/O_L} - \vec{\alpha}_{L/O_L}\times\vec{r}_{TMDI/P_L} - \vec{\omega}_{L/O_L}\times\left(\vec{\omega}_{L/O_L}\times\vec{r}_{TMDI/P_L}\right)$$
$$- 2\vec{\omega}_{L/O_L}\times\dot{\vec{r}}_{TMDI/P_L}. \quad (15)$$

in which the rotational acceleration terms are

$$\vec{\alpha}_{L/O_L}\times\vec{r}_{TMDI/P_L} = z_{TMDI_Z/P_L}\left[\ddot{\phi}_{L/O_L} \quad -\ddot{\theta}_{L/O_L} \quad 0\right]^T,$$
$$\vec{\omega}_{L/O_L}\times\left(\vec{\omega}_{L/O_L}\times\vec{r}_{TMDI/P_L}\right) = z_{TMDI_Z/P_L}\left[\dot{\theta}_{L/O_L}\dot{\psi}_{L/O_L} \quad \dot{\phi}_{L/O_L}\dot{\psi}_{L/O_L} \quad -(\dot{\theta}^2_{L/O_L}+\dot{\phi}^2_{L/O_L})\right]^T, \quad (16)$$
$$2\vec{\omega}_{L/O_L}\times\dot{\vec{r}}_{TMDI/P_L} = \dot{z}_{TMDI_Z/P_L}\left[2\dot{\phi}_{L/O_L} \quad -2\dot{\theta}_{L/O_L} \quad 0\right]^T.$$

Substituting Eqs. (14) and (16) into Eq. (15) and projecting along the local z direction gives the governing equation of motion for $TMDI_Z$ case as

$$\ddot{z}_{TMDI_Z/P_L} = \left[\frac{m_z(\dot{\theta}^2_{L/O_L}+\dot{\phi}^2_{L/O_L}) - k_z}{m_z + b_z}\right] z_{TMDI_Z/P_L} - \frac{c_z}{m_z + b_z}\dot{z}_{TMDI_Z/P_L} - \frac{m_z}{m_z + b_z}\ddot{z}_{P/O_L} - \frac{b_z}{m_z + b_z}(\ddot{z}_{P/O_L} - \ddot{z}_{T/O_L}) + \frac{m_z}{m_z + b_z}a_{G_Z/O_L} + \frac{1}{m_z + b_z}(F_{ext_Z} + F_{StopFrc_Z} + F_{Z_{PreLoad}}), \quad (17)$$

with constraint reactions

$$\begin{aligned} F_{X_{TMDI_Z/O_L}} &= m_z[-a_{G_X/O_L} + \ddot{x}_{P/O_L} + (\ddot{\phi}_{L/O_L} + \dot{\theta}_{L/O_L}\dot{\psi}_{L/O_L})z_{TMDI_Z/P_L} + 2\dot{\phi}_{L/O_L}\dot{z}_{TMDI_Z/P_L}], \\ F_{Y_{TMDI_Z/O_L}} &= m_z[-a_{G_Y/O_L} + \ddot{y}_{P/O_L} - (\ddot{\theta}_{L/O_L} - \dot{\phi}_{L/O_L}\dot{\psi}_{L/O_L})z_{TMDI_Z/P_L} - 2\dot{\theta}_{L/O_L}\dot{z}_{TMDI_Z/P_L}]. \end{aligned} \quad (18)$$

## 2.6 State-space representtion and OpenFAST Integration

To enable numerical integration within OpenFAST, the second-order equations of motion in Eqs. (8), (12) and (17) for the three independent TMDI directions of action need to be reformulated as a system of first-order ordinary differential equations. To this aim, a state vector $\vec{R}_{TMDI/P_L}$ is defined to collect the TMDI positions and velocities along their respective degrees of freedom relative to point P expressed in the local orientation $L$

$$\vec{R}_{TMDI/P_L} = \left[x_{TMDI_x/P_L} \;\; \dot{x}_{TMDI_x/P_L} \;\; y_{TMDI_y/P_L} \;\; \dot{y}_{TMDI_y/P_L} \;\; z_{TMDI_z/P_L} \;\; \dot{z}_{TMDI_z/P_L}\right]^T. \quad (19)$$

The input state equation can then be written as

$$\dot{\vec{R}}_{TMDI/P_L} = \boldsymbol{A}\vec{R}_{TMDI/P_L} + \boldsymbol{B}. \quad (20)$$

where the system matrix $\boldsymbol{A}$ is

$$\boldsymbol{A} = \begin{bmatrix} 0 & 1 & 0 & 0 & 0 & 0 \\ \frac{m_x(\dot{\phi}^2_{L/O_L}+\dot{\psi}^2_{L/O_L}) - k_x}{m_x + b_x} & -\left(\frac{c_x}{m_x + b_x}\right) & 0 & 0 & 0 & 0 \\ 0 & 0 & 0 & 1 & 0 & 0 \\ 0 & 0 & \frac{m_y(\dot{\theta}^2_{L/O_L}+\dot{\psi}^2_{L/O_L}) - k_y}{m_y + b_y} & -\left(\frac{c_y}{m_y + b_y}\right) & 0 & 0 \\ 0 & 0 & 0 & 0 & 0 & 1 \\ 0 & 0 & 0 & 0 & \frac{m_z(\dot{\theta}^2_{L/O_L}+\dot{\phi}^2_{L/O_L}) - k_z}{m_z + b_z} & -\left(\frac{c_z}{m_z + b_z}\right) \end{bmatrix}, \quad (21)$$

and the input vector $\boldsymbol{B}$ is

$$\boldsymbol{B} = \begin{bmatrix} 0 \\ \frac{-m_x\ddot{x}_{P/O_L} - b_x(\ddot{x}_{P/O_L} - \ddot{x}_{T/O_L})}{m_x + b_x} + \frac{m_x}{m_x + b_x}a_{G_x/O_L} + \frac{1}{m_x + b_x}(F_{ext_x} + F_{StopFrc_X}) \\ 0 \\ \frac{-m_y\ddot{y}_{P/O_L} - b_y(\ddot{y}_{P/O_L} - \ddot{y}_{T/O_L})}{m_y + b_y} + \frac{m_y}{m_y + b_y}a_{G_y/O_L} + \frac{1}{m_y + b_y}\left(F_{ext_y} + F_{StopFrc_Y}\right) \\ 0 \\ \frac{-m_z\ddot{z}_{P/O_L} - b_z(\ddot{z}_{P/O_L} - \ddot{z}_{T/O_L})}{m_z + b_z} + \frac{m_z}{m_z + b_z}a_{G_z/O_L} + \frac{1}{m_z + b_z}(F_{ext_z} + F_{StopFrc_Z} + F_{Z_{PreLoad}}) \end{bmatrix}. \quad (22)$$

The output vector $\vec{Y}_{P_G} = [\vec{F}_{P_G} \; \vec{M}_{P_G}]^T$ collects the resulting forces $\vec{F}_{P_G}$ and moments $\vec{M}_{P_G}$ acting on point P in global orientation, expressed as

$$\vec{F}_{P_G} = R_{L/G}^T \begin{bmatrix} k_x x_{TMDI_x/P_L} + c_x \dot{x}_{TMDI_x/P_L} - F_{ext_x} - F_{StopFrc_X} - F_{X_{TMDI_y/O_L}} - F_{X_{TMDI_z/O_L}} \\ k_y y_{TMDI_y/P_L} + c_y \dot{y}_{TMDI_y/P_L} - F_{ext_y} - F_{StopFrc_Y} - F_{Y_{TMDI_x/O_L}} - F_{Y_{TMDI_z/O_L}} \\ k_z z_{TMDI_z/P_L} + c_z \dot{z}_{TMDI_z/P_L} - F_{ext_z} - F_{StopFrc_Z} - F_{Z_{PreLoad}} - F_{Z_{TMDI_x/O_L}} - F_{Z_{TMDI_y/O_L}} \end{bmatrix},$$

$$\vec{M}_{P_G} = R_{L/G}^T \begin{bmatrix} -\left(F_{Z_{TMDI_y/O_L}}\right) y_{TMDI_y/P_L} + \left(F_{Y_{TMDI_z/O_L}}\right) z_{TMDI_z/P_L} \\ \left(F_{Z_{TMDI_x/O_L}}\right) x_{TMDI_x/P_L} - \left(F_{X_{TMDI_z/O_L}}\right) z_{TMDI_z/P_L} \\ -\left(F_{Y_{TMDI_x/O_L}}\right) x_{TMDI_x/P_L} + \left(F_{X_{TMDI_y/O_L}}\right) y_{TMDI_y/P_L} \end{bmatrix}, \tag{23}$$

respectively. Further, the total forces and moments acting on point T by the inerter are given as

$$\vec{F}_{T_G} = R_{L/G}^T \begin{bmatrix} b_x\left(\ddot{x}_{TMDI_x/P_L} + \ddot{x}_{P/O_L} - \ddot{x}_{T/O_L}\right) \\ b_y\left(\ddot{y}_{TMDI_y/P_L} + \ddot{y}_{P/O_L} - \ddot{y}_{T/O_L}\right) \\ b_z\left(\ddot{z}_{TMDI_z/P_L} + \ddot{z}_{P/O_L} - \ddot{z}_{T/O_L}\right) \end{bmatrix}, \qquad \vec{M}_{T_G} = R_{L/G}^T \begin{bmatrix} 0 \\ 0 \\ 0 \end{bmatrix}, \tag{24}$$

respectively.

Equations (19) to (24) are implemented in the StC sub-module. At each time step of simulation, the component kinematics $\ddot{\vec{r}}_{P/O_L}$, $\vec{\omega}_{L/O_L}$, $\dot{\vec{\omega}}_{L/O_L}$, and $\ddot{\vec{r}}_{T/O_L}$ are retrieved from the structural solver via the glue code interface. These kinematic quantities, together with the current TMDI state vector $\vec{R}_{TMDI/P_L}$, are used to evaluate the system matrix $\boldsymbol{A}$ and input vector $\boldsymbol{B}$ in Eqs. (21) and (22). The state evolution in Eq. (20) is then integrated numerically to advance the TMDI states to the next time-step. Subsequently, the output forces $[\vec{F}_{P_G}\ \vec{F}_{T_G}]$ and moments $[\vec{M}_{P_G}\ \vec{M}_{T_G}]$ are computed from Eqs. (23) and (24), and passed back to the structural modules through the glue code, thereby enabling fully coupled aero-hydro-servo-elastic, time-domain analysis of TMDI-equipped wind turbines as shown in Fig. 3.

Overall, the above implementation extends the registry file of StC module to accommodate six TMDI-specific parameters: three inertance values along x, y, and z axes and three spatial coordinates defining the inerter terminal point location. Correspondingly, six computational subroutines are modified to enable input parsing, parameter validation, dynamic mesh creation, state-derivative computation, and output load calculation. The implementation employs a dynamic two-node mesh architecture that creates the inerter terminal node only when TMDI functionality is activated, ensuring computational efficiency. Notably, backward compatibility with conventional TMD configurations is preserved that reduces TMDI equations to TMD equations when inertance values are set to zero. Full implementation details, including subroutine-level modifications, are detailed in Appendix A.

## 3. Numerical validation of TMDI implementation in OpenFAST

To demonstrate the accuracy and versatility of the achieved TMDI implementation in OpenFAST, illustrative numerical results are reported in this section, pertaining to TMDI placement at three different turbine components, namely turbine tower, blades, and floating substructure (Fig. 3). The focus is on comparing inerter force time-histories obtained from OpenFAST, against time-histories obtained from independent MATLAB® modeling and numerical integration. To this aim, the TMDI is modelled in MATLAB® environment as a mechanical network comprising a secondary mass connected to two nodes with externally prescribed kinematics through spring, damper, and inerter elements arranged as depicted in Fig. 1(b). The governing equation of motion for the TMDI network is written as

$$m\ddot{u}_m + c(\dot{u}_m - \dot{u}_1) + k(u_m - u_1) - b(\ddot{u}_2 - \ddot{u}_m) = 0, \tag{25}$$

where $m$, $b$, $k$, and $c$ denote the mass, inertance, stiffness and damping of TMDI, respectively, taken equal to the TMDI parameters reported in all the ensuing OpenFAST implementations. Further, in Eq. (25), $u_m$ represents the displacement of the TMDI mass, while $u_1$ and $u_2$ denote the prescribed displacements at terminals P and T of the network, respectively, in the modified OpenFAST StC sub-module. Displacement, velocity, and acceleration time histories at the TMDI attachment points P and T from OpenFAST are applied to the TMDI network terminals of the MATLAB model and inerter force time-histories from OpenFAST are compared with those derived by MATLAB. For each different TMDI implementation, a 500 s long time-domain simulation is run for external metocean excitations described below, and inerter force time-histories are compared for an initial 5 s-long transient time window, followed by a 15 s-long transition time window leading to steady-state response.

The first illustrative example considers a tower-top TMDI implementation to the IEA 15 MW offshore reference wind turbine supported by bottom-fixed monopile foundation detailed in [44]. The publicly available OpenFAST model of the considered IEA 15MW benchmark structure [45] is herein modified to accommodate two independent TMDIs at the tower-top. An OpenFAST executable compiled from the extended source code detailed in Section 2 and in Appendix A is used for all simulations conducted in this work. With reference to Fig. 3(a), one TMDI is acting in the fore-aft (FA) direction x and one in the side-side (SS) direction y. The TMDIs have identical properties: mass $m$ = 322.2 kg corresponding to 0.01% of the total mass of the IEA 15MW turbine including the monopile foundation equal to 3222 Mg [44]; inertance $b$=322.2 Mg, corresponding to 10% of the total structural mass; inerter span $h$= 28.9 m, corresponding to 20% of the tower height from the transition piece; and stiffness $k$ = 399.04 kN/m and damping coefficient $c$ = 42.225 kNs/m, derived from optimal TMDI tuning for the considered structure and TMDI inertial properties (i.e. mass and inertance) as detailed later in Section 4.2. For the OpenFAST simulation, metocean loading conditions for a power production design load combination (DLC) with extreme sea state, classified as DLC 1.6 per IEC 61400-3-1 guidelines [46], is used. Wind speed at hub height is taken equal to the turbine rated wind speed of 10.59 m/s [44], combined with extreme wave loading with 50-year return period consistent with the SEM-REV testing site, offshore France [47], as detailed in Appendix B. For numerical integration, the built-in fourth-order Adams-Bashforth-Moulton algorithm of OpenFAST and fourth-order Runge-Kutta algorithm of MATLAB have been used with a constant time-step of $5x10^{-4}$ s.

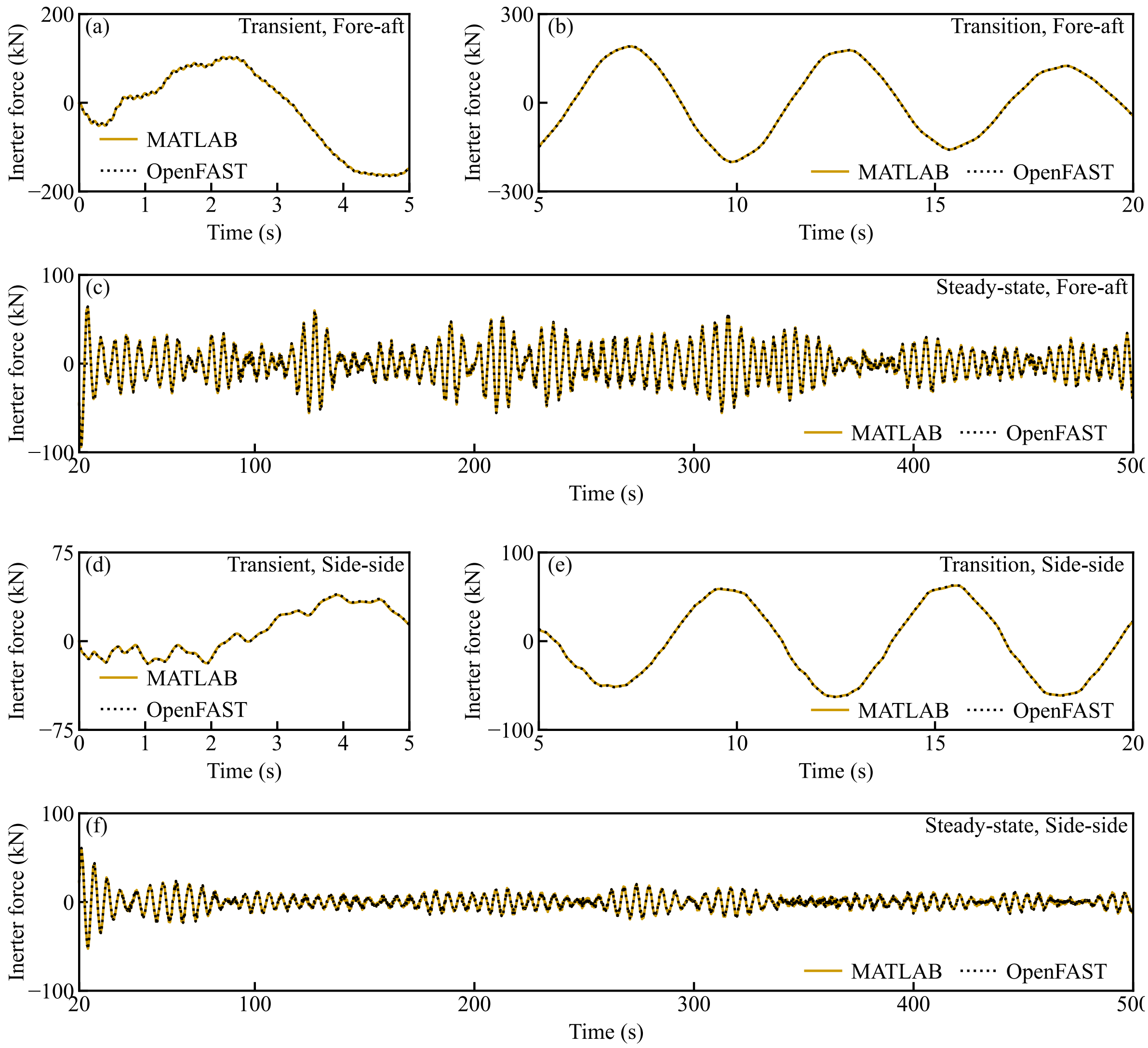

**Fig. 5** Inerter force time histories of the tower-top TMDI across transient, transition, and steady-state response phases, computed using OpenFAST and MATLAB for time step dt=0.0005s: (a) transient FA response, (b) transition FA response, (c) steady-state FA response, (d) transient SS response, (e) transition SS response, and (f) steady-state SS response.

Independently derived inerter force time-histories from OpenFAST and MATLAB are plotted in Fig. 5 for the full 500 s of simulation for the two TMDIs in the FA and SS directions. Excellent agreement is observed across all three different time windows considered (transient, transition, and steady state), showcasing the successful and accurate implementation of TMDI dynamics in OpenFAST for tower-top TMDI configuration. The potential of using OpenFAST to examine the capability of TMDI for mitigating OWT tower vibrations is exemplified in the following Section.

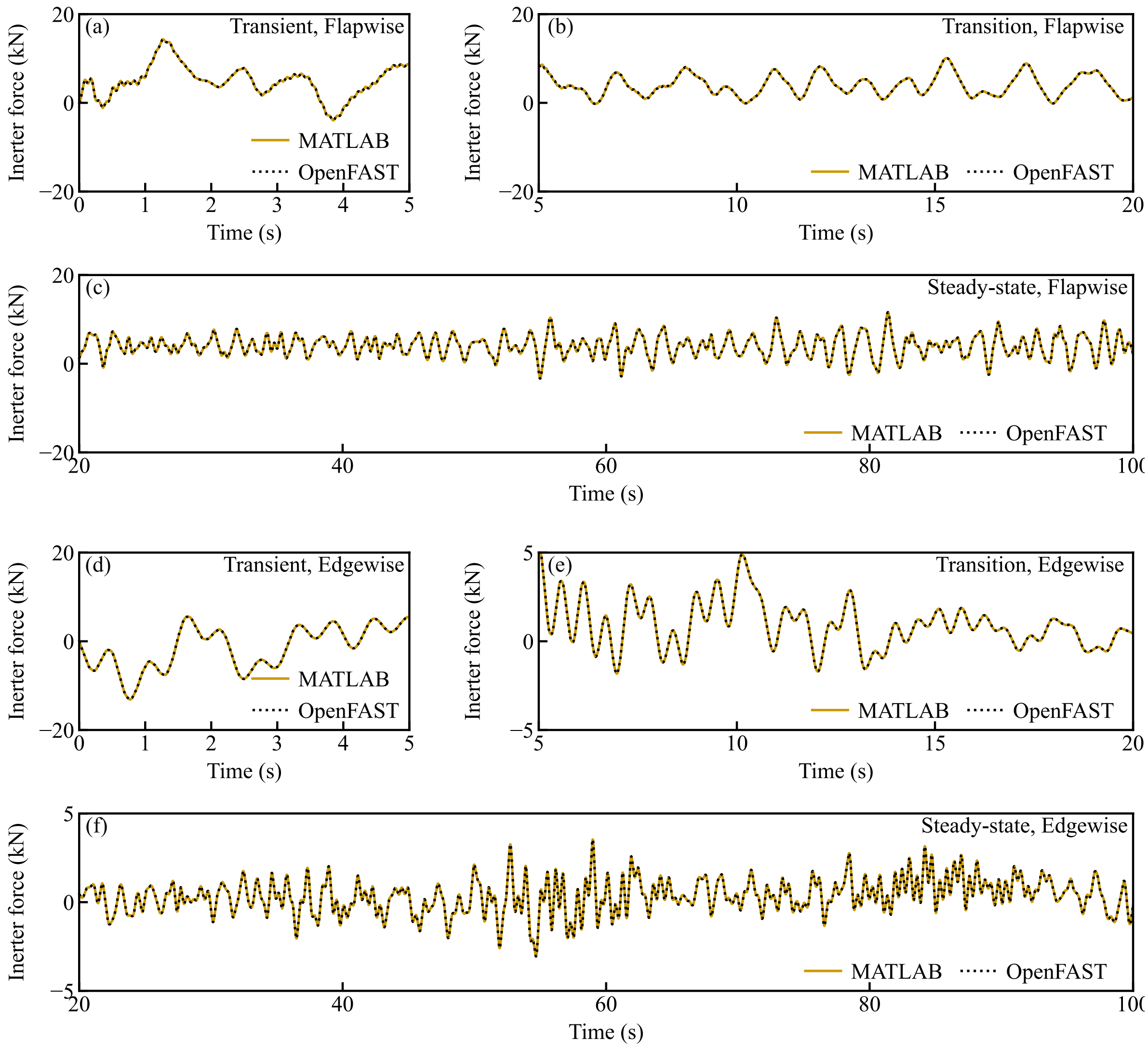

**Fig. 6** Inerter force time histories of the blade-mounted TMDI across transient, transition, and steady-state response phases, computed using MATLAB and OpenFAST for time step dt = 0.0005 sec: (a) transient flapwise response, (b) transition flapwise response, (c) steady-state flapwise response, (d) transient edgewise response, (e) transition edgewise response, and (f) steady-state edgewise response.

The second illustrative example considers an in-blade TMDI implementation to the same benchmark turbine IEA 15 MW. The OpenFAST model of the turbine is modified to include two independent TMDIs in each blade, one in the local flapwise x direction and one in the local edgewise y direction (Fig. 3 (c)) with attachment points P at the root of the blades. The TMDIs have identical properties: mass $m = 6.84$ kg corresponding to 0.01% of the mass of each blade equal to 68.4 Mg [44]; inertance $b =$ 6.84 Mg, corresponding to 10% of the blade mass; inerter span $h = 28.9$ m, as taken in the previous example; and stiffness $k$ = 75.18 kN/m and damping coefficient $c$ = 5.65 kNs/m, taken as arbitrary but reasonable values for the considered application. The same metocean loading conditions are used for OpenFAST simulation as in the previous example with the same numerical integration algorithms and time step in OpenFAST and MATLAB. Inerter force time-histories from OpenFAST and MATLAB are plotted in Fig. 6 for the first 100 s of simulation for two TMDIs in the flapwise and edgewise directions, demonstrating the same level of excellent agreement as in the previous example (Fig. 5) during the transient, transition, and steady state response. These results showcase the feasibility of using OpenFAST to assess the potential of TMDI for suppressing blade oscillations in large-scale wind turbines instead of proprietary wind turbine models previously considered in the literature [36].

As a final example, a TMDI OpenFAST implementation in a sub-structure is illustrated using the OC4 DeepCwind semi-submersible floating platform supporting the reference NREL 5 MW offshore wind turbine [48], as a benchmark model structure. For this case, OpenFAST model files of the benchmark floating OWT [49] are modified to include one TMDI inside each of the three columns of the semi-submersible platform along the vertical (heave) direction z (Fig. 3(b)). The TMDIs are attached to the inner column walls at 12 m below the top of the column (i.e. d=12 m in Fig. 3(b)), with spring, dashpot, and inerter in parallel (vertical) orientation (i.e. d = 12 m and $h = 0$ in Fig. 3(b)). Therefore, this example examines the limiting case of a single-point TMDI in Fig. 1(a) with $\ddot{u}_2$ replaced by $\ddot{u}_1$ in Eq. (25). The three TMDIs have identical properties: mass $m$ = 137.23 kg corresponding to 0.01% of the total mass of the wind turbine plus the floating platform equal to 1372.3 Mg; inertance $b$ = 137.23 Mg, corresponding to 10% of the total mass; stiffness $k$ = 23 kN/m and damping coefficient $c$ = 1.72 kNs/m, taken as arbitrary but reasonable values for the considered application. For metocean loading conditions, a wind speed of 12 m/s at 90 m hub height is assumed, while waves are characterized by the JONSWAP spectrum (see Appendix B, with significant wave height $H_s$ = 1.79 m, peak period $T_p$ = 10 s, and peak shape parameter $\gamma = 1$. Further, a smaller time step of $10^{-4}$ s is adopted for numerical integration in OpenFAST and MATLAB, compared to the previous two examples, necessitated to accurately capture the impact-like TMDI dynamics at the start of the simulation under gravitational self-weight loading along the heave direction. The latter is manifested by the large magnitude inerter force acting in the heave direction at the start of the transient response shown in Fig. 7. Overall, juxtaposed inerter force time-histories from OpenFAST and MATLAB shown in Fig. 7 are in excellent agreement as in the previous examples, showcasing the applicability of using OpenFAST to investigate motion control potential in floating substructures and platforms.

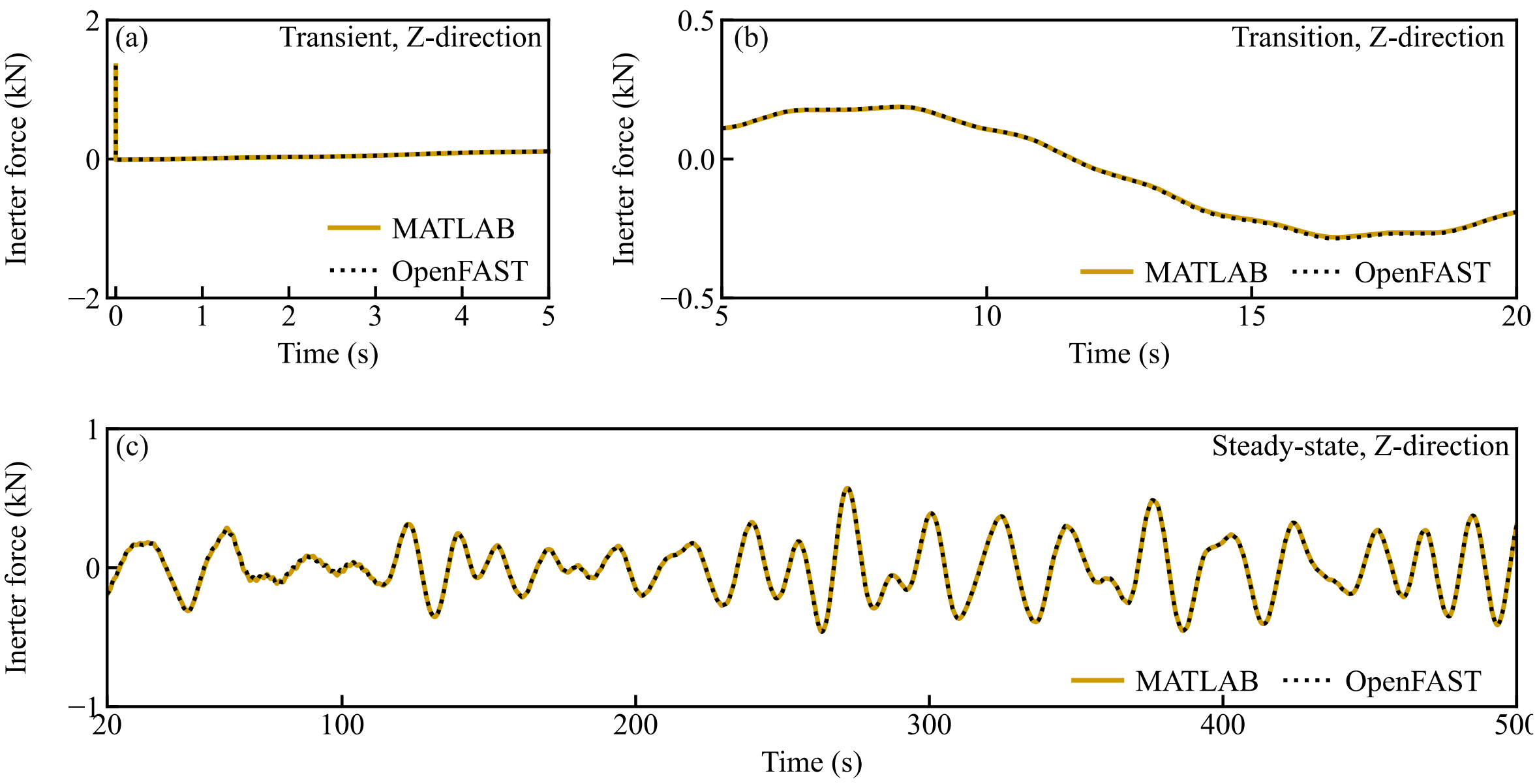

**Fig. 7** Inerter force time histories of the semi-submersible-platform TMDI in the heave (vertical) direction, computed using MATLAB and OpenFAST for time step dt = 0.0001 sec: (a) transient response, (b) transition response, and (c) steady-state response..

## 4. Illustrative OWT assessment in OpenFAST fitted with optimally tuned TMDI at tower-top

Having established the applicability and accuracy of implementing the TMDI in OpenFAST, this Section provides an illustrative case-study to quantify the potential of TMDI versus the standard TMD for mitigating the tower oscillations in large-scale monopile OWTs under different turbine operational conditions. For this purpose, the IEA 15MW reference turbine on monopile foundation is used [44] under power production and idling conditions, equipped with bi-directional tower-top lightweight TMDIs of different inertance (see also Fig. 3(a)), as well as a standard bi-directional TMD. The main aim is to assess the benefits of increasing inertance to the performance of TMDI-equipped wind turbines, following previous studies (e.g. [30–32,40]), but using time-domain OpenFAST coupled aero-hydro-servo-elastic simulations to do so, for the first time in the literature. This is supported by frequency/excitation-agnostic optimal TMD(I) tuning which safeguards fairness in the performance comparison of the different passive vibration absorbers considered in parametric investigation. The presentation begins with the description and finite element modeling of the benchmark structure used to obtain the dominant mode shapes and natural frequencies of the tower-plus-monopile system which are required for the optimal TMDI tuning.

### 4.1 Description and finite element modelling of the IEA 15MW reference offshore wind turbine

Specified in [44], the IEA 15MW reference OWT is a three-bladed upwind horizontal-axis wind turbine supported by a monopile foundation with a fixed base at the mudline as shown in Fig. 8. The rotor has a diameter of 240 m, with the hub positioned at 150 m above mean sea level (MSL). The tower extends approximately 130 m from the MSL to the tower top, with an outer diameter tapering from 10 m at the base to 6.5 m at the nacelle interface. The monopile has a total length of 75 m and a constant outer diameter of 10 m, of which 45 m is embedded in the seabed and the remaining 30 m extends through the water column.

The publicly available OpenFAST model of the IEA 15MW includes natural frequencies and mode shapes of only the OWT tower (i.e. above the transition piece) in the input files of the ElastoDyn module (see Fig. 2). However, it cannot furnish the coupled tower-plus-monopile modal properties needed for an accurate TMD(I) tuning. To this end, a stick-leg finite element model (FEM) of the considered OWT is developed in OpenSeesPy finite element software [50,51], whose accuracy is validated against OpenFAST natural frequencies and mode shapes of the tower. In the developed FEM, the tower-plus-monopile structure is discretized using 120 two-node displacement-based beam-column elements with 5 integration points, as depicted in Fig. 8. The outer diameter is constant at 10 m from the mudline to the transition piece, followed by linear tapering to the tower top, whilst wall thickness reduces stepwise from the mudline upwards, following the IEA 15 MW specifications in [42]. Further, the OWT transition piece and the rotor-nacelle assembly (RNA) are modelled as lumped masses at the monopile–tower intersection node and the tower-top node, respectively. At the tower top, a rigid offset is used to account for the eccentricity of the hub height with respect to the tower top, as shown in the inlet of Fig. 9, following the actual geometry of the OWT detailed in [44]. The adopted steel material properties in the FEM are: modulus of elasticity equal to 210 GPa; Poisson's ratio equal to 0.3; and mass density equal to 7.85 $Mg/m^3$ which is increased by 7% to account for the mass contribution of connections and welding [44].

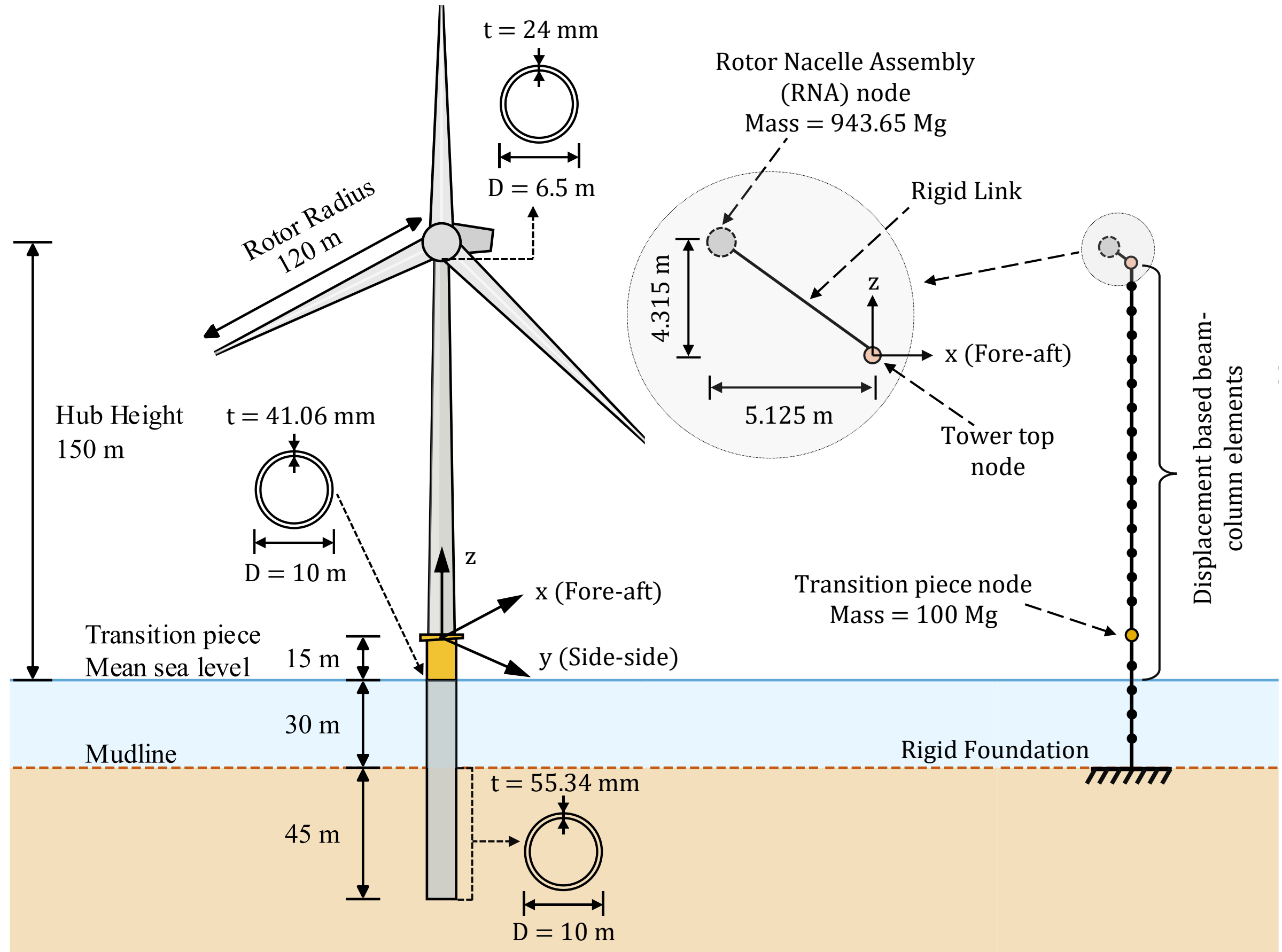

**Fig. 8** Schematic representation of the simplified finite element model of the IEA-15 MW reference offshore wind turbine in OpenSeesPy.

The accuracy of the developed FEM is verified by comparing tower-only natural frequencies and mode shapes, vis-à-vis those coming from the ElastoDyn module input files of the OpenFAST turbine model. The FEM-based properties are obtained by undertaking standard modal analysis in OpenSeesPy, without accounting for geometric nonlinearity. Fig. 9 plots the first two tower mode shapes in FA and SS directions obtained from OpenSeesPy and from OpenFAST, while Table 1 reports the natural frequencies of the above mode shapes alongside modal assurance criterion (MAC) values computed as

$$\mathrm{MAC} = \frac{|\,\varphi_{OF}^{T}\varphi_{OS}\,|^{2}}{(\varphi_{OF}^{T}\varphi_{OF})(\varphi_{OS}^{T}\varphi_{OS})}, \tag{26}$$

where $\varphi_{OF}$ and $\varphi_{OS}$ are mode shape vectors from OpenFAST and OpenSeesPy models, respectively, and the superscript "$T$" denotes matrix transposition. The MAC ratio in Eq. (26) is a common metric to quantify the degree of similarity between modes. The data in Fig. 19 and in Table 1 evidence excellent agreement in modal properties estimated by the FEM developed, thus validating its applicability for providing tower-plus-monopile modal properties required for optimal TMD(I) tuning.

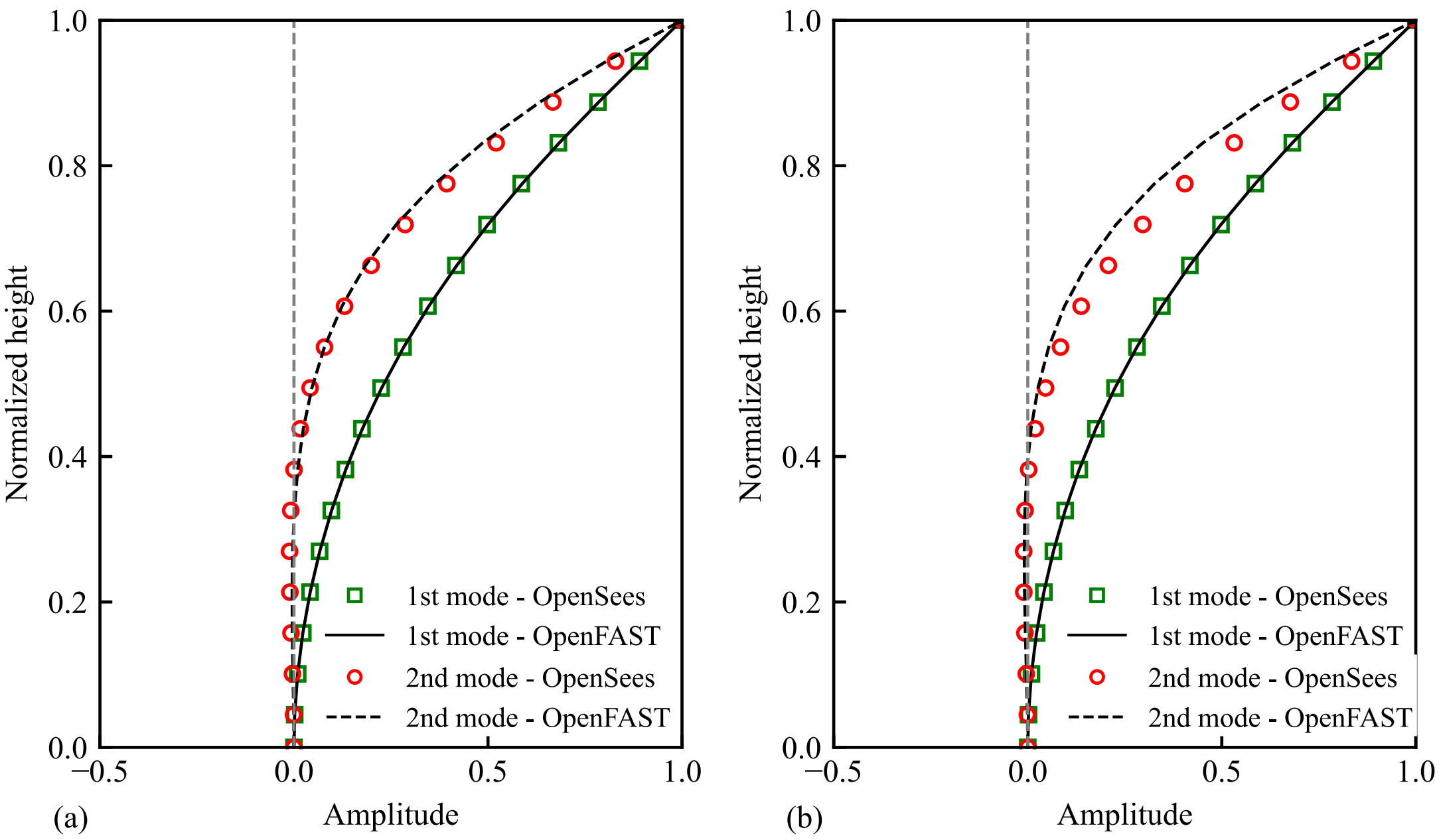

**Fig. 9** Comparison of normalized tower-only mode shapes of IEA 15MW obtained from OpenSeesPy and OpenFAST. (a) Fore-aft modes (b) Side-side modes.

**Table 1** Comparison of first two natural frequencies in FA and SS directions obtained from OpenFAST and OpenSeesPy for the IEA 15MW turbine [44].

| Mode | OpenFAST (Hz) | OpenSeesPy (Hz) | MAC |
|---|---|---|---|
| $1^{st}$ side-side | 0.175 | 0.179 | 1.000 |
| $1^{st}$ fore-aft | 0.179 | 0.179 | 1.000 |
| $2^{nd}$ side-side | 1.276 | 1.279 | 0.993 |
| $2^{nd}$ fore-aft | 1.225 | 1.257 | 0.999 |

## 4.2 Optimal TMD(I) tuning

In this study, all the considered TMD(I)s are tuned to suppress the tower-top displacement of the IEA 15MW turbine, using one, the first/fundamental, mode shape of the tower-plus-monopile structure assuming fixed support at mudline. In this setting, the tower-plus-monopile structure is represented as a generalized single-degree-of-freedom system, following [52,53], with coordinate corresponding to the tower-top displacement, $\bar{u}(t)$, critical damping ratio $\xi_1$, and natural frequency $\omega_1$ corresponding to the fundamental mode shape $\varphi_1(z)$ as delineated in Fig. 10(a). This mode shape and natural frequency are obtained from the FEM of the tower-plus-monopile structure developed in OpenSeesPy and validated against OpenFAST (Fig. 10 and Table 1). Importantly, the first mode along the FA and SS directions turns out to be practically the same based on modal analyses applied to the FEM. Therefore, identical modal properties and thus tuning of TMD(I)s along FA and SS directions are taken.

With the addition of a TMD(I) to the OpenSeesPy-based generalized single-degree-of-freedom system, a two-degree-of-freedom system is obtained with $\bar{u}(t)$ and $u_d(t)$ coordinates (degrees-of-freedom) as shown in Fig. 10(c). The tower-top displacement frequency response function (FRF) of the latter system, normalized by the static tower-top displacement under distributed force excitation, is given as [52]

$$H(\omega) = \frac{(\nu^2 - g^2 + j\,2g\nu\xi_d)(\mu+\beta)}{[1 - g^2(1+\beta\varphi_1(\chi)^2) + (\nu^2+2jg\nu\xi_d)(\mu+\beta) + 2jg\xi_1](\nu^2 - g^2 + j\,2g\nu\xi_d)(\mu+\beta) - [(\nu^2+2jg\nu\xi_d)(\mu+\beta) - g^2\beta\varphi_1(\chi)]^2} \quad (27)$$

where $j = \sqrt{-1}$, $g = \omega/\omega_1$ is a normalized frequency, $\chi$ is the height of the inerter attachment point T onto the turbine measured from the mudline (Fig. 10), and the mass ratio $\mu$, the inertance ratio $\beta$, the TMDI frequency ratio $\nu$, and the TMDI damping ratio $\xi_d$ are defined as

$$\mu = \frac{m}{M}, \beta = \frac{b}{M}, \nu = \frac{\sqrt{k/(m+b)}}{\omega_1}, \xi_d = \frac{c}{2\sqrt{(m+b)k}}, \quad (28)$$

respectively. In the last expressions, $M = 3222$ Mg is the total mass of the IEA 15MW wind turbine.

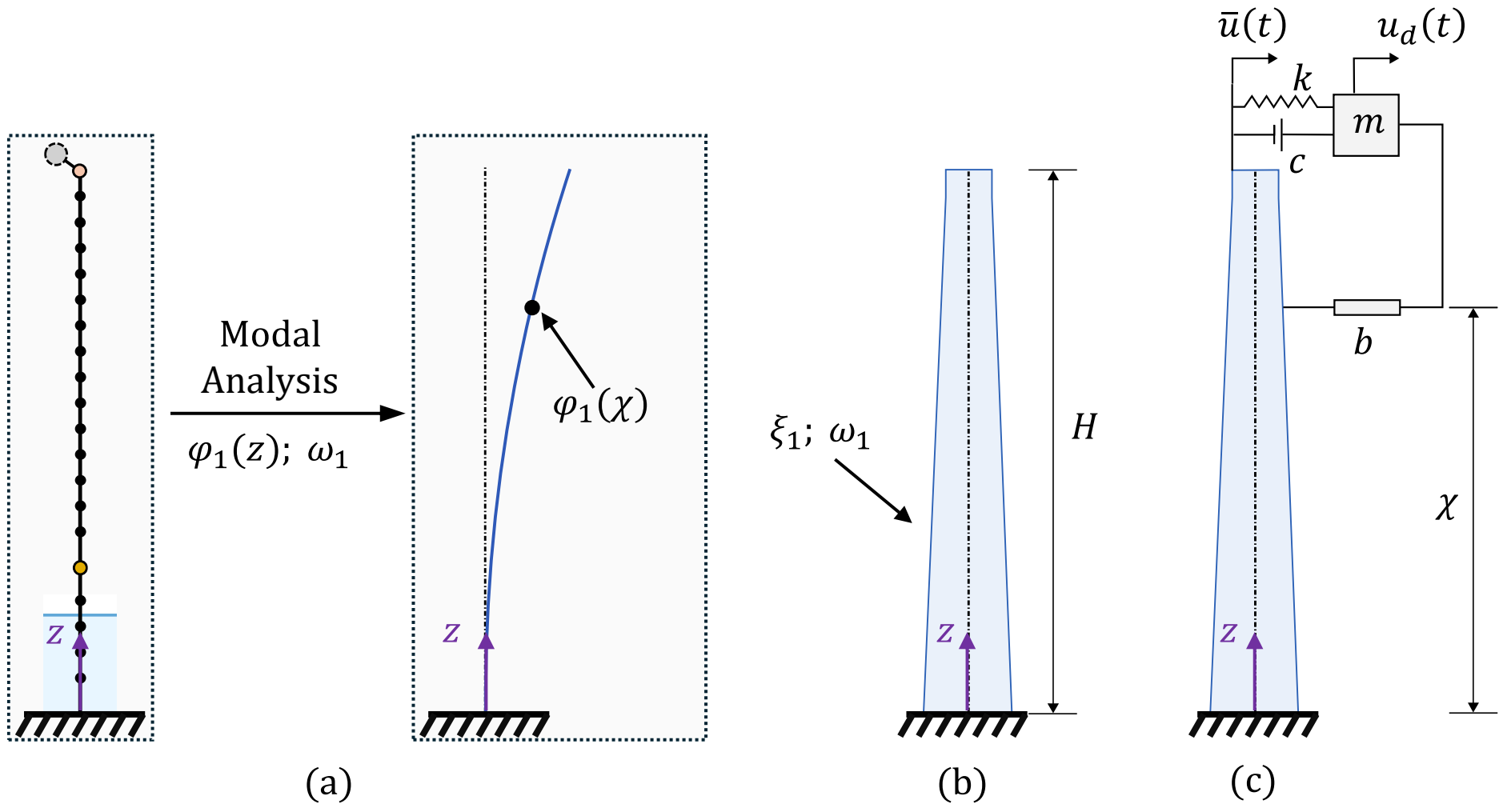

**Fig. 10** Two degree of freedom system modelling framework for TMDI-equipped IEA-15 MW OWT (a) OpenSeesPy finite element model of IEA-15 MW and extracted normalized mode shape (b) generalized SDOF representation of the uncontrolled primary structure (c) TMDI-equipped generalized SDOF system.

To this end, TMD(I) tuning is formulated as an $H_2$-style optimization problem, aiming to determine optimal TMDI frequency ratio, $\nu$, and damping ratio, $\xi_d$ that minimizes the area under the squared magnitude of the nondimensional FRF in Eq. (27), given mass ratio $\mu$, inertance ratio $\beta$, and inerter connectivity point location $\chi$. Mathematically, this optimization problem can be written as

$$A_{min} = \min_{\mathbf{y}_1} \left\{ \int_0^{\omega_c} |\, H(\omega)\,|^2\, d\omega \right\} \text{given } \mathbf{y}_2 \text{ and subjected to } \mathbf{y}_1^{\min} \leq \mathbf{y}_1 \leq \mathbf{y}_1^{\max}, \quad (29)$$

in which $\omega_c$ is the cut-off frequency above which $H(\omega)$ attains negligible values, $\mathbf{y}_1 = [\nu, \xi]^T$ is the vector of primary TMDI design variables with lower and upper bounds specified by $\mathbf{y}_1^{\min}$, and $\mathbf{y}_1^{\max}$, respectively, and $\mathbf{y}_2 = [\mu, \beta, \chi]^T$ is the vector of secondary (pre-specified) TMDI design variables. Notably, the above optimal tuning approach can be used for TMD tuning by setting $\beta = b = 0$ in Eqs. (27)- (29).

The optimization problem in Eq. (29) is solved numerically using a gradient-free pattern search algorithm implemented in Python with search bounds $\mathbf{y}_1^{\min} = [0.0{,}0.0]^T$ and $\mathbf{y}_1^{\max} = [2.0, 2.0]^T$. The algorithm employs an iterative refinement strategy wherein the search domain is progressively reduced around the provisional optimal values of $\nu$ and $\xi$ from each iteration. Convergence is achieved when the relative change in objective function between successive iterations is below $10^{-2}$. For illustration, Fig. 11 shows the convergence behavior of the pattern search algorithm for optimal tuning of a TMDI with $\mu = 0.01\%$, inertance ratio $\beta = 10\%$, inerter connectivity $\chi = 0.8H$, attached to the IEA 15MW primary structure approximated the first mode of the OpenSeesPy FEM in Section 4.1 and assuming critical damping ratio $\xi_1$ =1%. The latter figure also exemplifies the strong convexity of the optimization problem in Eq. (29) encountered for all different TMD(I)s discussed next.

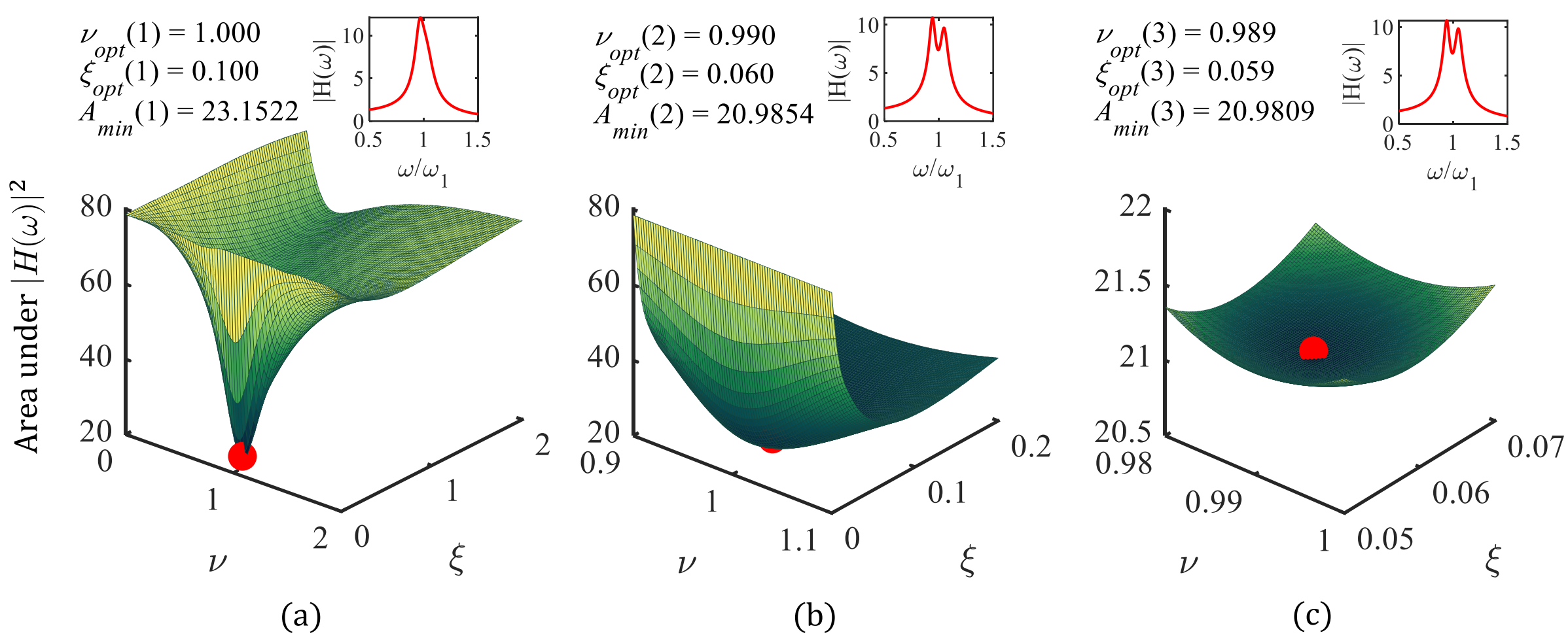

**Fig. 11** Convergence behavior of the pattern search algorithm for the $H_2$ optimization problem in Eq. 47. System parameters: $\mu = 0.01\%$, $\beta = 10\%$, $\chi = 0.8H$, $\xi_1 = 1\%$ and fore-aft mode shape $\varphi_1(z)$ obtained from OpenSeesPy FE model of IEA-15-240-RWT. The red dot in each panel marks the optimal TMDI parameters ($\nu_{opt}$, $\xi_{opt}$) and the corresponding minimum objective function value $A_{min}$. Panels show: (a) first iteration, (b) second iteration, and (c) third iteration.

For a parametric investigation of the importance of the inertance to the performance of TMDI-equipped IEA 15MW OWT with $\xi_1$ =1%, five TMDIs with different inertance ratios ranging from $\beta = 5\%$ to 40%, fixed inerter connectivity $\chi = 0.8H$, and fixed relatively small mass ratio of $\mu = 0.01\%$ are considered, whereby the TMDI mass is included to model the self-weight of the TMDI devices and connections. Furthermore, a TMD with substantial mass, corresponding to $\mu = 1\%$ which is typical of current real-life OWT installations, is also included in the investigation to demonstrate potential benefits of TMDI over current TMD implementations. Table 2 lists optimal TMD(I) tuning parameters obtained by solving Eq. (29) as detailed above, while Fig. 12 plots the magnitude of the nondimensional top-tower displacement FRF, $|H(\omega)|$, in Eq. (27) for all the optimally tuned TMD(I)s. It is seen that for sufficiently large inertance ratio, $\beta \geq 10\%$, the lightweight TMDI with mass ratio $\mu = 0.01\%$ outperforms the 100 times heavier TMD in terms of suppressing top-tower displacement FRF at ±10 % the resonant frequency of the uncontrolled structure. However, note that these FRFs are excitation-independent and are derived from the simplified FEM of the benchmark OWT. Therefore, the next sub-section assesses the performance of the IEA 15MW equipped with the thus-obtained optimal TMD(I)s using OpenFAST time-domain simulations for different turbine operational and metocean conditions.

**Table 2** Optimal stiffness and damping parameters for TMD and TMDIs considered in the case study.

| Inertial absorber | Mass Ratio $\mu$ (%) | Inertance Ratio $\beta$ (%) | Optimal stiffness $k_{opt}$ (kN/m) | Optimal damping $c_{opt}$ (kNs/m) |
|---|---|---|---|---|
| TMD | 1 | 0 | 40.1 | 3.59 |
| TMDI | 0.01 | 5 | 201.86 | 15.09 |
| TMDI | 0.01 | 10 | 399.04 | 42.33 |
| TMDI | 0.01 | 20 | 781.63 | 117.84 |
| TMDI | 0.01 | 30 | 1148.44 | 212.86 |
| TMDI | 0.01 | 40 | 1499.68 | 322.58 |

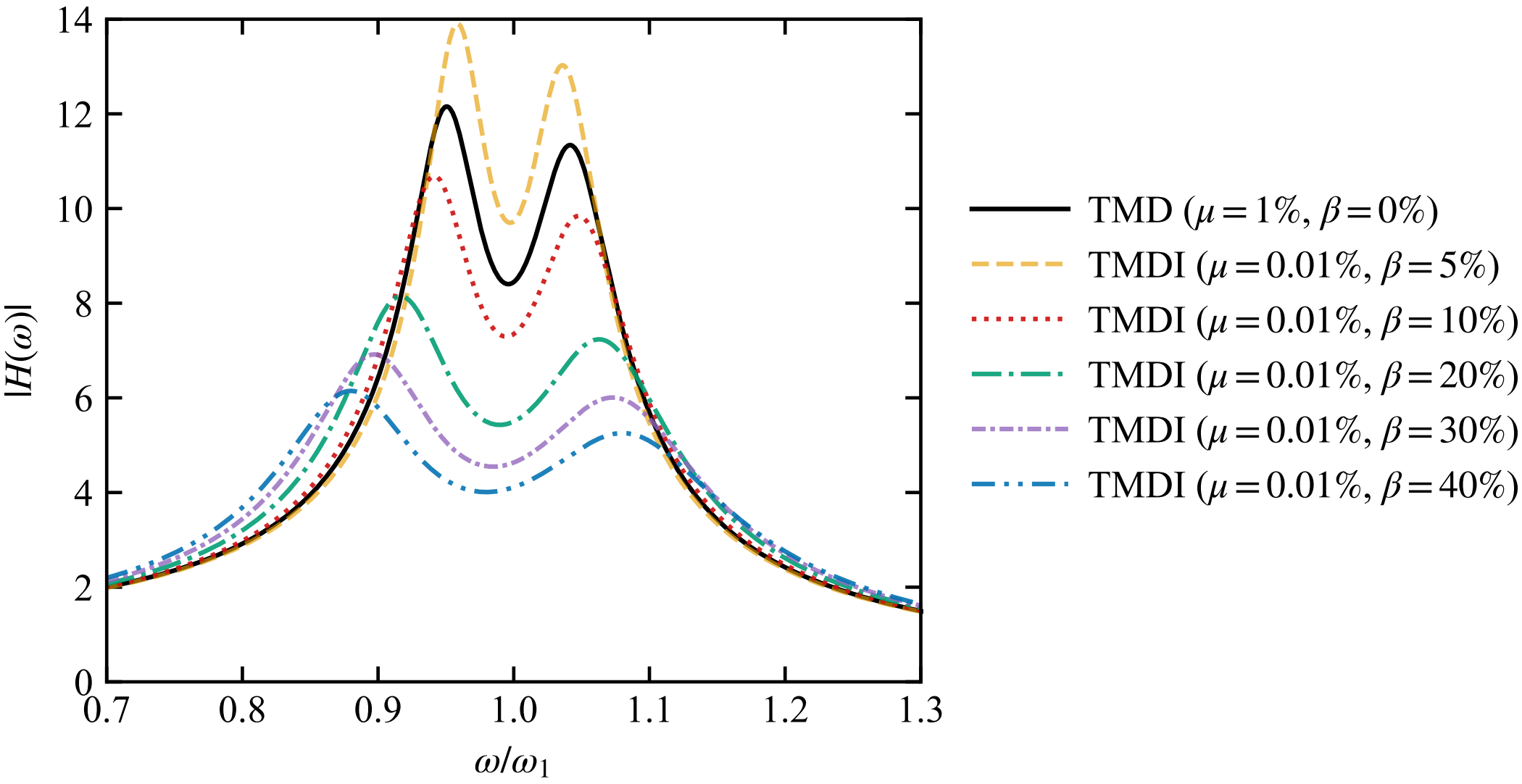

**Fig. 12** Magnitude of the nondimensional FRF, $|H(\omega)|$, for optimally tuned TMD and TMDIs under white noise excitation. The TMD configuration ( $\mu = 1\%$ ) is compared with TMDI configurations having $\mu = 0.01\%$ and inertance ratios $\beta = 5 - 40\%$. The primary structure corresponds to a generalized SDOF representation of IEA 15MW OWT with fore-aft mode shape $\varphi_1(z)$ and inherent critical damping ratio $\xi_1 = 1\%$.

### 4.3 Assessment of IEA 15MW equipped with optimized TMDI

Using the optimal TMD(I) tuning parameters of Table 2, performance assessment of IEA-15 MW fitted with one bi-directional TMD and five different bi-directional TMDIs is carried out by considering OpenFAST time-domain simulations spanning 700 s, with the initial 100 s of transient response discarded. Metocean loading conditions for two different DLCs are used consistent with the SEM-REV testing site, offshore France [47], as detailed in Appendix B: one power production DLC with normal sea state and rated wind speed at hub height, classified as DLC 1.2 per IEC 61400-3-1 guidelines [46]], relevant for fatigue limit state assessment; and one DLC for an idling (parked) wind turbine state under extreme storm conditions with 50-year return period, classified as DLC 6.1 per IEC 61400-3-1 guidelines [46], relevant for ultimate limit state assessment.

First, wind turbine performance is characterized in the frequency-domain by computing the power spectral density (PSD) of tower-top displacement response from OpenFAST simulation output, plotted and presented in Fig. 13. It is observed that spectral characteristics of the uncontrolled OWT depend strongly on both loading direction and operational state. For normal power production conditions (DLC 1.2), a pronounced resonant peak arises at the first natural frequency (~ 0.18 Hz) predominantly in the SS direction, whilst the FA response remains largely suppressed owing to aerodynamic damping provided by rotor rotation (Fig. 13(a), 13(b)). In contrast, under idling turbine and storm conditions (DLC 6.1), a prominent resonant peak emerges in the FA direction in the absence of aerodynamic damping, whilst a less pronounced resonant peak persists in the SS direction, as evidenced in Figs. 13(c) and 13(d). Importantly, TMD(I)s significantly modify these spectral characteristics. Under normal power production conditions, corresponding to DLC 1.2, the pronounced peak response in SS direction is effectively suppressed by both the conventional TMD and the TMDIs, whereas only a modest reduction is observed in the response along the FA direction. For DLC 6.1, TMD(I) effectiveness is more pronounced in the FA direction, whilst the improvement in the SS direction remains comparatively modest. Overall, the TMDIs with $\beta \geq 10\%$ consistently outperform the classical TMD with $\mu = 1\%$ across all examined cases, with performance improving monotonically as the inertance increases. This trend is consistent with FRFs in Fig. 12, proving that the relevance and efficiency of the FEM modelling and the low-order two-degree-of-freedom system used for TMD(I) optimal tuning presented in previous sections. Further, it directly demonstrates the beneficial effect of inertance combined with the TMDI configuration which outperforms a TMD with 100 times larger physical mass.

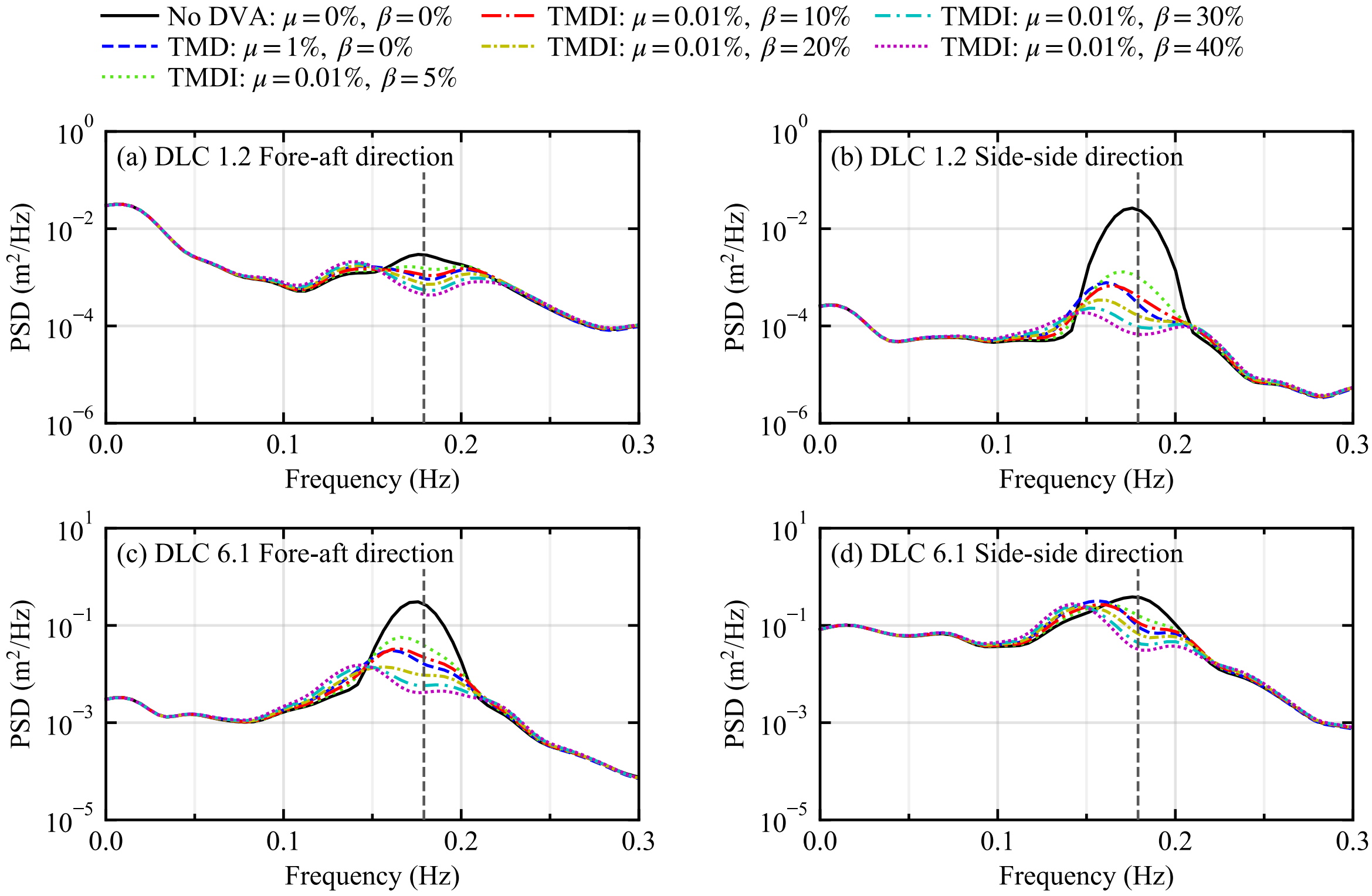

**Fig. 13** Power spectral density (PSD) of the fluctuating component of tower-top displacement response for the IEA-15 MW OWT considering all TMD(I)s of Table 2: (a) FA direction for DLC 1.2, (b) SS direction for DLC 1.2, (c) FA direction for DLC 6.1, and (d) SS direction for DLC 6.1. The vertical dashed line indicates the first natural frequency of the OWT.

Next, standard deviation and peak tower-top FA and SS displacement values are reported from OpenFAST response time-history data for DLC 1.2 in Fig. 14 and DLC 6.1 in Fig. 15 for all the TMD(I) cases considered, as well as for the uncontrolled (i.e. with no TMD(I)) wind turbine. On the same figures, percentage reductions with respect to the uncontrolled wind turbine are reported for the TMDI(I) cases. For normal power production conditions (DLC 1.2- Fig. 14), TMDI(I) effectiveness varies strongly with direction. The FA response shows minimal sensitivity to vibration control, with standard deviation reductions of 0.8% – 3.6% for $\beta = 5 - 40\%$ and peak displacement reductions below 2%. On the antipode, the SS response reduces substantially with TMD(I) addition, with TMD achieving 62% standard deviation reduction and a 17.5% peak reduction, while TMDIs attain reductions up to 67.3% for standard deviation and up to 18.9% for peak response. Notably, a TMDI with only 5% of inertance achieves similar reduction with the 100 times heavier TMD. The directional differences reflect the dominant influence of aerodynamic damping during power production in the FA direction which makes TMDI(s) ineffective and, ultimately, unnecessary. In contrast, the inherent wind turbine damping in the SS direction is much lower, making the response susceptible to resonance effects which can be efficiently controlled by the TMD(I).

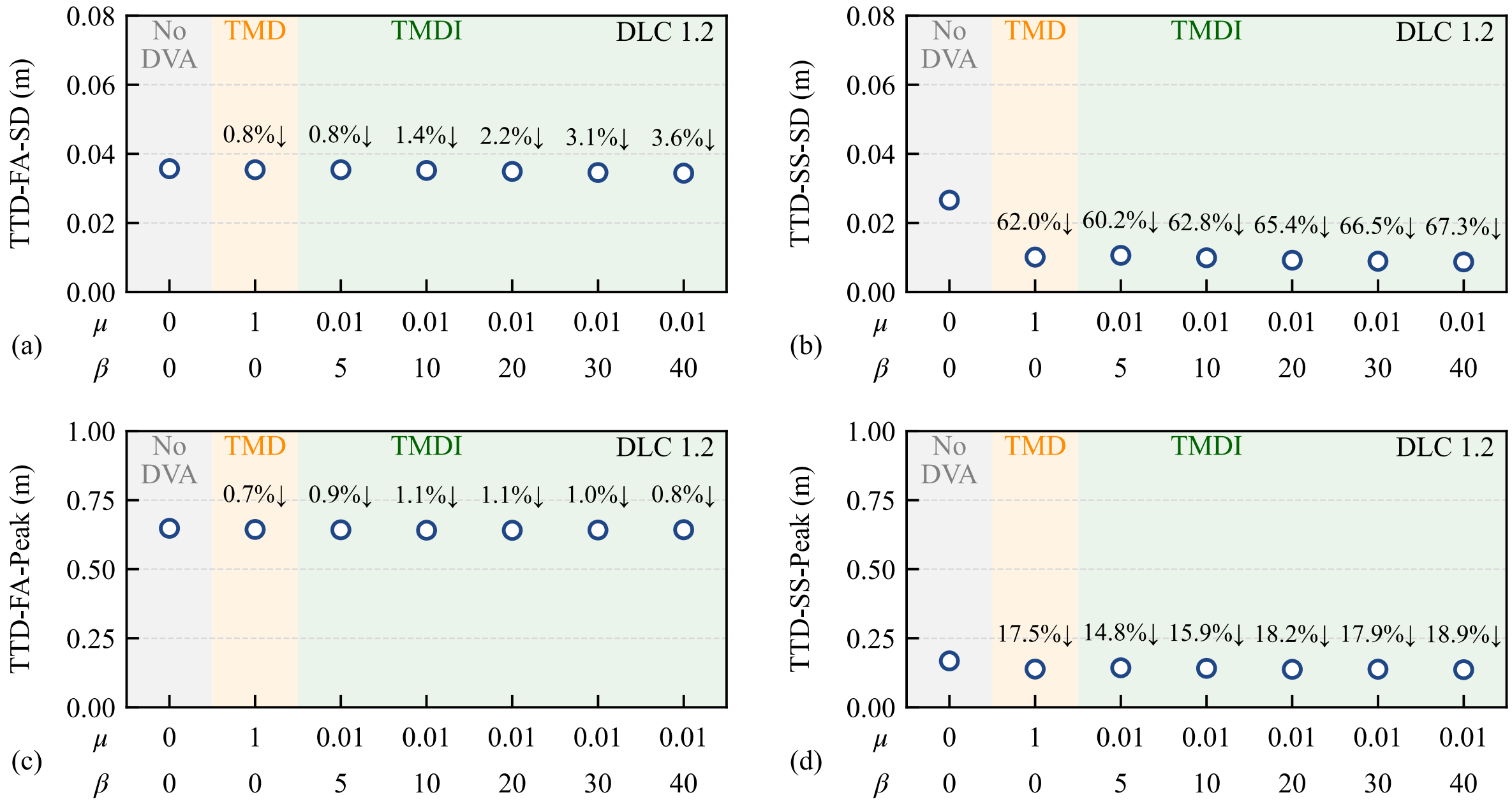

**Fig. 14** Tower-top displacement response for all TMD(I)s of Table 2 under DLC 1.2 operational conditions: standard deviation in the top row, with (a) FA direction and (b) SS direction; and peak response in the bottom row, with (c) FA direction and (d) SS direction. Percentage values with downward arrows indicate reduction relative to uncontrolled response.

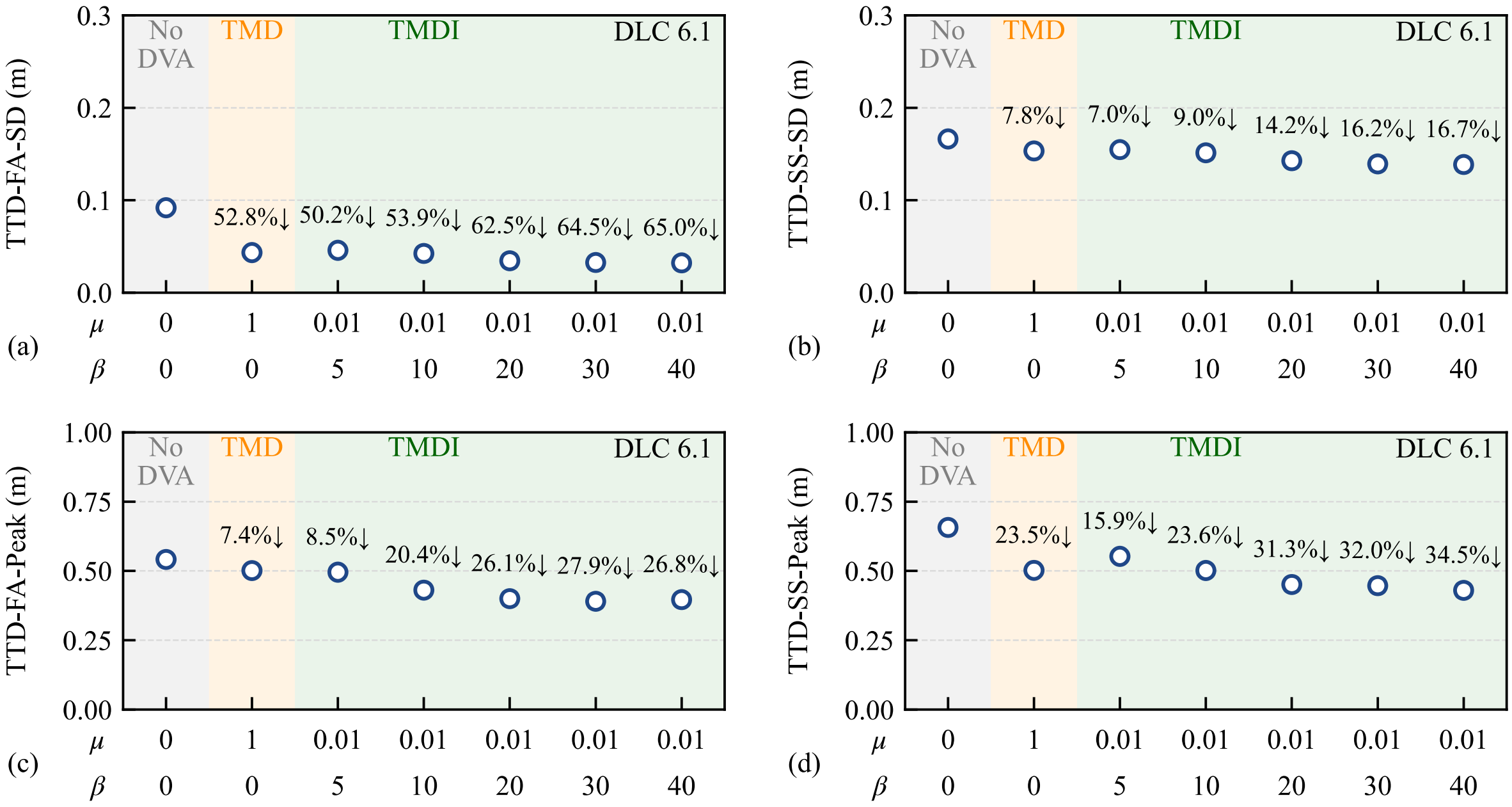

**Fig. 15** Tower-top displacement response for all TMD(I)s of Table 2 under DLC 6.1 idling conditions: standard deviation in the top row, with (a) FA direction and (b) SS direction; and peak response in the bottom row, with (c) FA direction and (d) SS direction. Percentage values with downward arrows indicate reduction relative to uncontrolled response.

Nevertheless, under extreme environmental conditions and turbine idling (DLC 6.1-Fig. 15), TMD(I) motion control effectiveness is appreciable in both FA and SS directions due to the absent of direction-dependent aerodynamic damping. For the FA response, the TMDI with $\beta = 40\%$ achieve 65.0% standard deviation reduction and 27.9% peak reduction compared to 52.8% and 7.4% for TMD, respectively, and similar performance reductions are exhibited in SS direction. Further, TMDIs with inertance $\beta \geq 10\%$ consistently exceed the performance of the 100-fold heavier TMD.

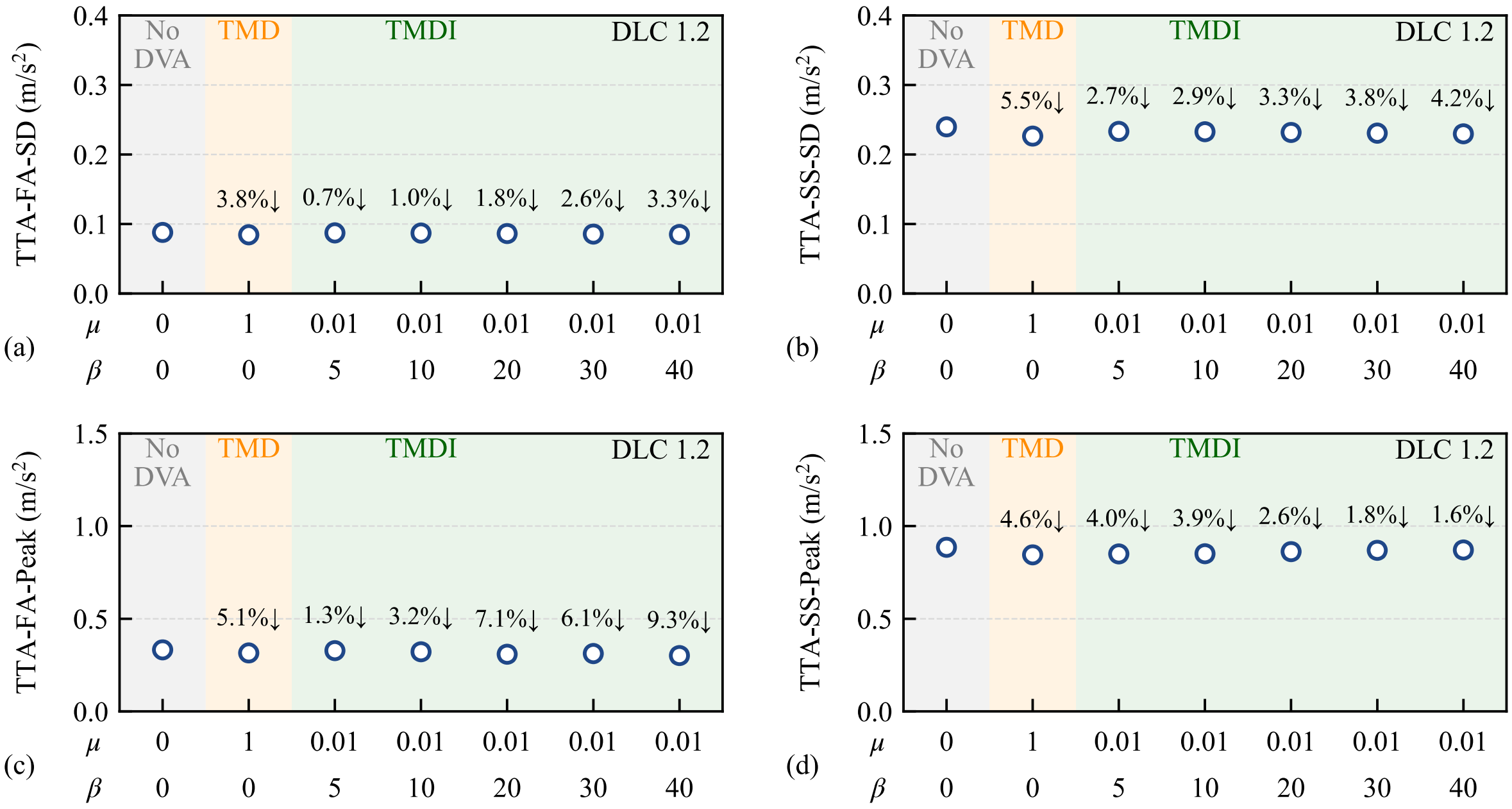

**Fig. 16** Tower-top acceleration response for all TMD(I)s of Table 2 under DLC 1.2 operational conditions: standard deviation in the top row, with (a) FA direction and (b) SS direction; and peak response in the bottom row, with (c) FA direction and (d) SS direction. Percentage values with downward arrows indicate reduction relative to uncontrolled response.

Further to tower-top displacement response, standard deviation and peak tower-top acceleration response are examined for both DLCs in Figs. 16 and 17, as acceleration governs serviceability limit state design and operational availability. This is because excessive RNA accelerations compromise the functionality of sensitive equipment at the nacelle and may trigger unscheduled shutdowns [54]. For normal operational conditions (DLC 1.2- Fig. 16), TMD(I)s provide modest acceleration mitigation in both FA and SS directions, with classical TMD outperforming TMDIs, except for the case of peak accelerations in FA direction. Nevertheless, for the extreme loading scenario with idling turbine (DLC 6.1-Fig. 17), the TMD(I)s become quite effective in top-tower acceleration mitigation, and their effectiveness depends largely on the inertance property. For instance, the TMDI with $\beta = 40\%$ achieves reductions in the standard deviation acceleration of 60.8% and 46.2% in FA and SS directions, respectively, which are 30% and 70% higher than the corresponding reductions that TMD achieves along the two directions. Further, TMD(I)s are also effective in mitigating peak tower-top accelerations, in both FA and SS directions. In particular, the TMDI with $\beta = 10\%$ achieves a 27.5% reduction in peak acceleration along the SS direction, markedly outperforming the corresponding 12.2% reduction attained by the TMD. Conversely, in the FA direction, TMDI peak acceleration reduction performance is broadly comparable to that of the TMD, with the highest TMDI reduction (46.4%, achieved at $\beta = 10\%$) only marginally exceeding the TMD's 45.8% reduction. In this respect, collectively with the case of tower-top displacement, it is concluded that TMDI with at least 10% of inertance provides significant supplemental damping to mitigate tower kinematics when aerodynamic damping is absent due to turbine idling under extreme metocean conditions, thus directly contributing to wind turbine survivability.

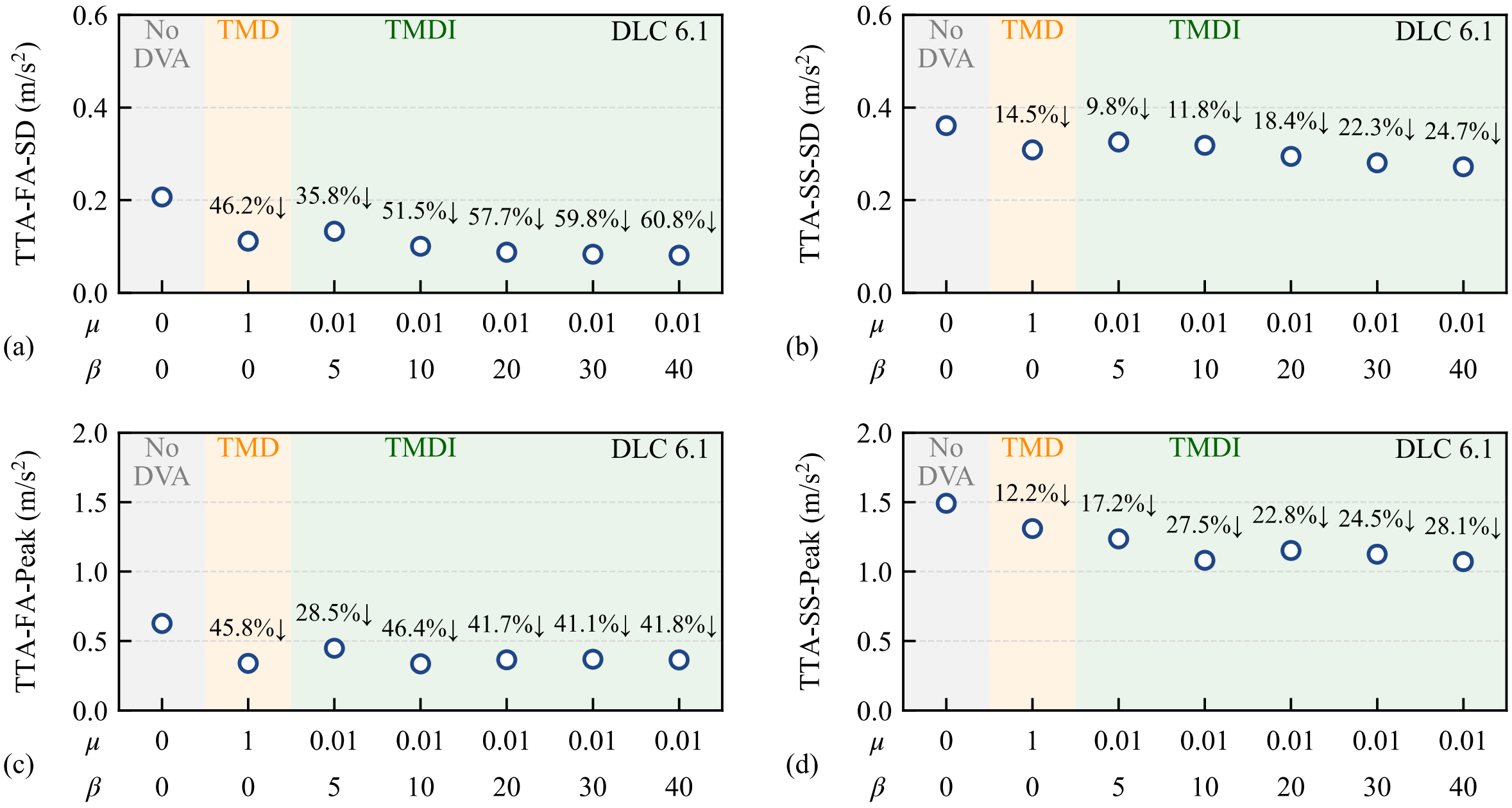

**Fig. 17** Tower-top acceleration response for all TMD(I)s of Table 2 under DLC 6.1 idling conditions: standard deviation in the top row, with (a) FA direction and (b) SS direction; and peak response in the bottom row, with (c) FA direction and (d) SS direction. Percentage values with downward/upward arrows indicate reduction/increment relative to uncontrolled response.

Beyond tower-top kinematics mitigation, which by-and-large improves fatigue and ultimate stresses developing along the turbine tower and monopile substructure, examining the TMD(I) peak stroke, defined as peak secondary mass displacement relative to the tower-top displacement, is also important. This is because the damper stroke relates to the practical implementation of TMD(I) given space limitations within the tower, as well as the TMD(I) costs (e.g. higher stroke results in higher viscous damper costs). Given that the IEA 15MW tower apex has a diameter of 6.5 m, a peak permissible stroke of 3 m is herein assumed, leaving a clearance of 0.25 m from the tower sidewalls. In practice, hard stoppers are used to limit the TMD stroke within permissible values, which are modelled in the StC module of OpenFAST by means of stiff double-sided impact viscoelastic elements with a pre-specified gap from the TMD secondary mass set equal to the permissible stroke [55,56]. In this respect, large TMD stroke, beyond permissible values, result in pounding/collision to the stoppers which compromise the motion control capability of TMDs and induce high instantaneous impact forces to the tower. Fig. 18 reports the peak strokes for all the considered

TMD(I)s and for both DLCs. For normal operational conditions (DLC 1.2) peak stroke demands are low and well within the permissible limit. However, TMD(I) stroke demands increase significantly under extreme loading conditions and idling turbine (DLC 6.1). In fact, the TMD exceeds the permissible stroke of 3 m in SS direction under DLC 6.1 (note that peak stroke reaches 3.54 m for TMD in SS direction obtained from OpenFAST simulation with no stoppers in the modelling) and collides to the stoppers as indicated by a dashed line in Fig. 18(d). Notably, such excessive TMD stroke, exceeding the available tower-top space, have not been previously reported in the literature, which has hitherto focused on relatively smaller-capacity wind turbines (e.g., NREL 5MW) with substantially lower stroke demands than the IEA 15MW. Nevertheless, the underlying trend of TMDI stroke reduction persists at this larger scale: TMDI strokes are appreciably lower than the TMD and are monotonically decreasing with inertance, which is in alignment with previously reported trends in the literature pertaining to the NREL 5MW wind turbine [30,40] and, more generally, for any slender cantilevered structure under lateral dynamic loading [52]. In this respect, it is herein confirmed using OpenFAST time-domain simulation data that the addition of an inerter to the TMD is beneficial not only for suppressing tower response, but also for reducing secondary mass stroke.

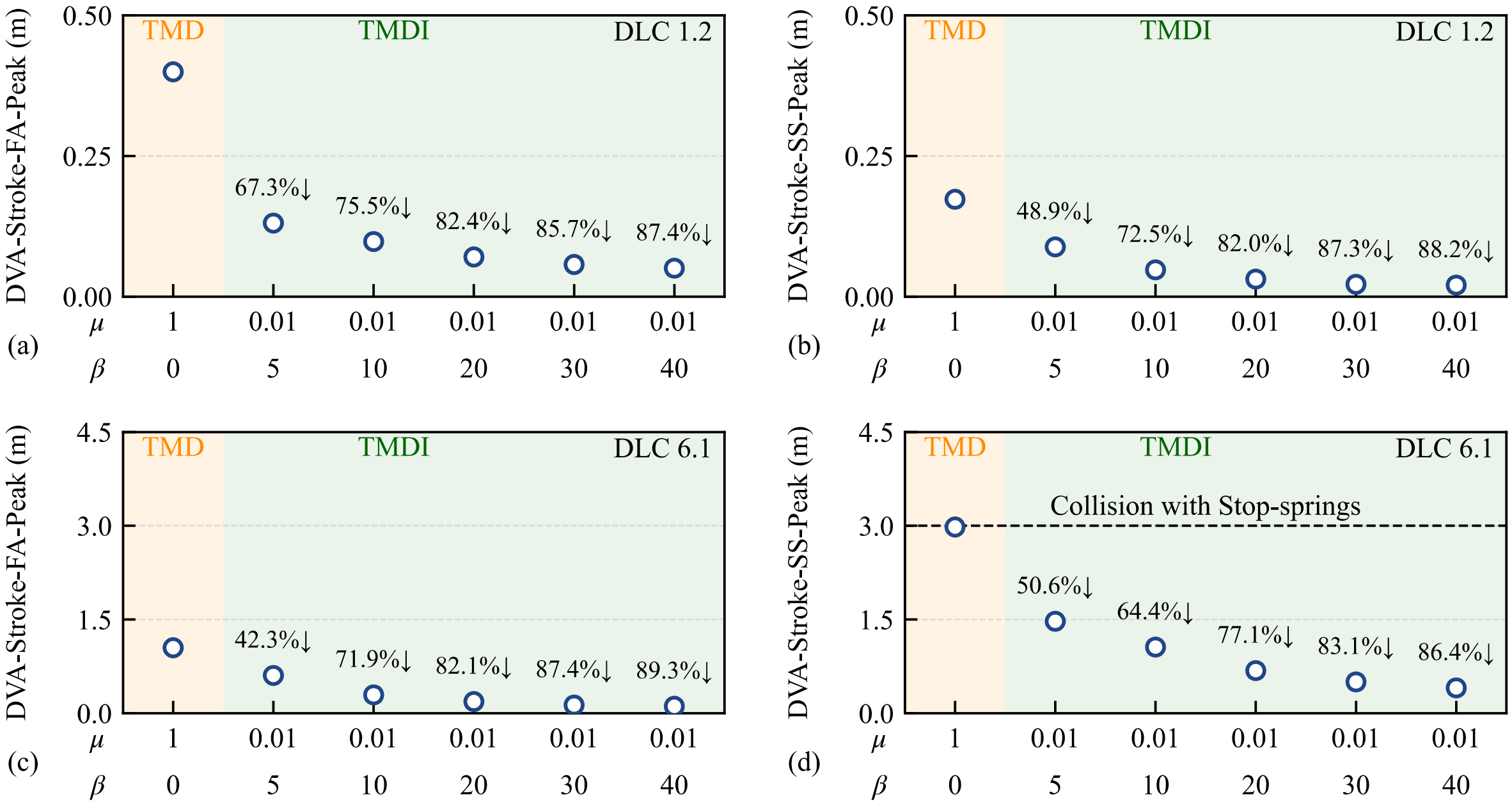

**Fig. 18** Peak values of secondary mass stroke for all TMD(I)s of Table 2 under DLC 1.2 operational conditions: (a) FA direction, (b) SS direction; and DLC 6.1 idling conditions: (c) FA direction, (d) SS direction. Percentage values with downward arrows indicate reduction relative to TMD. Dashed line in (d) represents stroke limit of 3 m at which collision of TMD mass takes place against stop-springs.

Nevertheless, previous studies considering the TMDI for mitigating the response of slender cantilevered structures under severe dynamic excitations showed that reductions in structural vibrations and damper stroke with increase inertance come at the cost of increased forces developing in the inerter device which need to be accommodated by the primary structure [57,58]. In fact, these inerter force demands have not been examined in the TMDI literature specific to wind turbines, which has primarily focused on vibration and stroke mitigation performance. In this respect, as a final assessment step in this illustrative case-study, it is deemed essential to examine peak inerter forces, which are reported in Fig. 19. It is found that inerter forces increase monotonically with inertance. Further, inerter force demands under extreme metocean conditions (DLC 6.1) are substantially higher than those developing under normal turbine operation (DLC 1.2) in both FA and SS directions. Still, even for the worst-case scenario of DLC 6.1, inertance ratio of 40%, and SS direction, the peak inerter force does not exceed 650 kN of force which is well within expected forces to develop in large-scale structures equipped with energy dissipation devices under severe dynamic environmental loading.

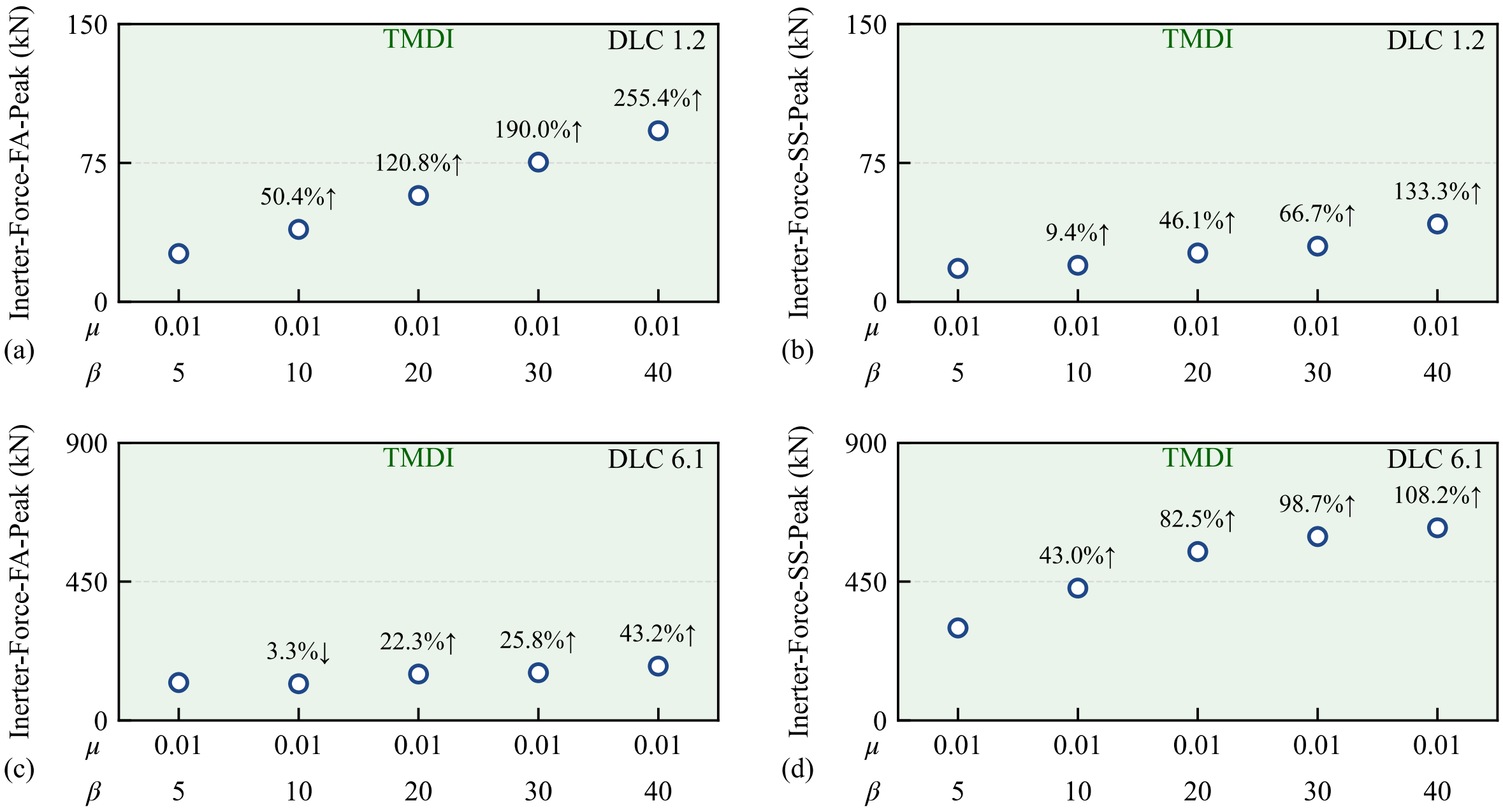

**Fig. 19** Peak values of Inerter force for all TMD(I)s of Table 2 under DLC 1.2 operational conditions: (a) FA direction, (b) SS direction; and DLC 6.1 idling conditions: (c) FA direction, (d) SS direction. Percentage values with downward/upward arrows indicate reduction/increment relative to TMDI ($\beta = 5\%$). TMDI: $\mu = 0.01\%$, $\beta = 5 - 40\%$.

Collectively, Figs. 14-19 confirm previous studies on TMDI in that higher inertance achieve superior vibration suppression, reduced weight requirements, and substantial stroke reduction compared to the TMD, at the cost of high inerter force demands which increase the device cost and require diligence at the inerter-tower connection. Therefore, the TMDI design challenge lies not only in inerter device capacity, but also in integrating these connections within the hollow cylindrical tower geometry while maintaining structural integrity under combined environmental and mechanical actions. Such challenges can only be effectively addressed through high-fidelity modelling at wind turbine design and assessment stages for accurate prediction of the TMDI-equipped OWT performance and demands (see also [59]). To this end, the herein extension of OpenFAST, addresses directly such requirements through high-fidelity time-domain simulations for assessing the dynamic response of TMDI-equipped OWTs.

## 5. Concluding Remarks

Tuned mass damper inerters (TMDIs) have been implemented and validated within OpenFAST's aero-hydro-servo-elastic simulation environment for the first time, enabling performance-driven design of lightweight vibration absorbers for wind turbines. A rigorous multibody dynamics formulation was developed and integrated into the structural control module, supporting arbitrary inerter connection locations across tower, nacelle, substructure, and blade components with independent multi-directional configurations.

The implementation was verified against independent MATLAB numerical models for tower-based, substructure-based, and blade-based TMDI configurations. It is found that excellent agreement is achieved across time windows considered thereby establishing the implementation's accuracy for TMDI design within OpenFAST's coupled simulation environment.

To demonstrate the practical applicability of TMDIs in OpenFAST, performance assessment was conducted for the IEA-15 MW offshore wind turbine equipped with a tower-top bi-directional TMDI under combined wind-wave loading. Optimal TMDI parameters were determined through $H_2$-norm minimization of tower-top displacement using a simplified two-degree-of-freedom model, with modal properties obtained from an independent OpenSeesPy structural model validated against OpenFAST. It is demonstrated that TMDIs with secondary mass ratio of 0.01% and inertance ratio of 10% achieve vibration suppression performance matching or exceeding that of a conventional TMD with 1% mass ratio across two critical design load cases under realistic coupled aero-hydro-servo-elastic loading: power production at rated conditions (DLC 1.2) and idling under storm conditions (DLC 6.1). These results confirm the mass amplification effect inherent to TMDI configurations, in which a physical mass 100 times smaller than the TMD mass generates equivalent control forces through inerter-enabled inertia enhancement. Furthermore, it is observed that TMDIs impart supplemental damping at frequencies exceeding 1 Hz beyond the fundamental mode, indicating potential for broadband vibration mitigation across multiple structural modes.

Importantly, the OpenFAST simulations demonstrate that TMDI stroke requirements are reduced by approximately 90% compared to conventional TMD under extreme loading conditions. This significant reduction enables practical integration within spatially constrained tower-top geometries while simultaneously reducing the required secondary mass by two orders of magnitude. In this regard, the developed capability facilitates performance evaluation of TMDIs as practical alternatives to conventional TMDs in applications where space limitations and weight constraints present significant design challenges.

Nevertheless, TMDI implementation entails a fundamental trade-off between vibration suppression effectiveness and inerter force demands, which increase monotonically with inertance. Higher inertance values, while enhancing performance, necessitate special structural design provisions at inerter connection points to accommodate the exerted forces locally. In this respect, site-specific optimization of TMDI parameters considering both performance objectives and structural force constraints is critical, which is directly facilitated by the developed implementation.

The integration of TMDI functionality within OpenFAST enables parametric investigation of optimal inerter connection locations and inertance values across different wind turbine configurations, spanning land-based and offshore systems with various support structures and turbine ratings. Furthermore, the open-source nature of the implementation promotes facilitates algorithmic development within the wind energy research community, thereby paving the way towards resilient renewable energy infrastructure.

## Appendix A: Source Code Modifications for TMDI Implementation in OpenFAST

Integration of the TMDI within OpenFAST is achieved through systematic modifications to the Structural Control (StC) module. The implementation strategy involves extensions to both the data structure definitions and the computational subroutines, as detailed in the following. Importantly, the modifications preserve backward compatibility with existing TMD configurations while enabling the enhanced capabilities afforded by TMDI.

### Registry File Modifications

The StC module architecture in OpenFAST is defined through a registry file (*StrucCtrl_Registry.txt*), which establishes the interface between user-specified input parameters and internal computational data structures. To accommodate TMDI functionality, three categories of modifications are introduced to this registry file.

First, the description of the degree-of-freedom mode parameter is updated to reflect the expanded translational capabilities of the module. Specifically, mode 1 is redefined to explicitly denote three independent translational degrees of freedom (X, Y, Z), each of which may optionally employ a TMDI configuration. This modification is applied consistently across both the input file data section and the internal parameters section of the registry.

Second, six additional input fields are introduced to specify TMDI-specific parameters. Three fields define the inertance values along the X, Y, and Z axes, measured in units of mass (kg), while the remaining three fields specify the spatial coordinates defining the location of the inerter terminal point T relative to the primary attachment point P, expressed in the local component orientation and measured in meters. These input fields are utilized exclusively when the translational degree-of-freedom mode is active.

Third, six corresponding internal parameters are added to store the TMDI configuration data during runtime simulations. These parameters comprise three inertance values and three spatial coordinates, which are mapped from the input file data during module initialization through the parameter-setting subroutine. Following these registry modifications, OpenFAST is recompiled to establish the foundational data architecture required for subsequent computational implementations.

**Computational Subroutine Modifications**

The computational logic of the StC module is implemented in the *StrucCtrl.f90* source file, which contains the subroutines governing structural control algorithms. Six subroutines are modified to accommodate TMDI functionality, organized into five functional categories: input file parsing, parameter mapping and validation, mesh initialization, state-derivative computation, and output load computation. The overall modular structure of the StC module remains unchanged, ensuring maintainability and extensibility.

Input file parsing: The input file parsing subroutine (*StC_ParseInputFileInfo*) is modified to read the six TMDI-specific parameters from the structural control input file. Three parsing operations retrieve the spatial coordinates of the inerter terminal point, while the remaining three extract the inertance values along the respective axes. These fields are parsed unconditionally but may be assigned zero values in input files when inerter functionality is not required, thereby maintaining backward compatibility with conventional TMD configurations.

Parameter mapping and validation: The parameter-setting subroutine (*StC_SetParameters*) is modified to map the input file data into internal parameter structures. This mapping is performed conditionally, transferring input values to parameter structures only when the appropriate translational mode is active. When alternative control modes are specified, TMDI-specific parameters are explicitly set to zero, ensuring data structure isolation and appropriate memory allocation. Furthermore, the parameter validation subroutine (*StC_ValidatePrimaryData*) is extended to enforce physically meaningful TMDI configurations. For each active translational degree of freedom, validation logic verifies that inertance values are non-negative and that mass and stiffness values are strictly positive. Violation of these constraints triggers a fatal error, preventing numerically ill-conditioned simulations.

Mesh initialization and node creation: The initialization subroutine (*StC_Init*) is modified to implement the two-node mesh architecture required for TMDI operation. In conventional TMD configurations, the StC module employs a single-node mesh representing the attachment point P. For TMDI configurations, a second node representing the inerter terminal point T is dynamically created when at least one non-zero inertance value is specified. This dynamic node allocation ensures that the computational overhead associated with the additional mesh node is incurred only when TMDI functionality is required.

The mesh creation logic is generalized to accommodate variable node counts, enabling dynamic allocation of either single-node or two-node mesh configurations. The mesh supports kinematic quantities including translational displacement, orientation, translational velocity, rotational velocity, translational acceleration, and rotational acceleration at each node. When a two-node mesh is instantiated, the global position of node T is computed by transforming the user-specified offset vector from local to global coordinates using the initial reference orientation matrix and superposing the result onto the position of node P. Using *MeshPositionNode* subroutine, Node T is positioned at this computed location with identical orientation to the reference frame, and a point element is constructed to enable application of forces and moments.

State-derivative computation: The continuous-time state-derivative subroutine (*StC_CalcContStateDeriv*) is modified to implement the TMDI acceleration formulations derived in Section 2. The acceleration vector at node T is first transformed from global to local coordinates through multiplication with the direction cosine matrix. Subsequently, state derivatives are computed using the governing equations presented in Eqs. 33-35. Notably, when all inertance values are set to zero, the formulation reduces automatically to the classical TMD equations. In this manner, a unified computational framework supports both conventional TMD and TMDI configurations, with behavior determined solely by the specified inertance values.

Output load computation: The output computation subroutine (*StC_CalcOutput*) is modified to evaluate forces and moments applied to the host component at points P and T. Local acceleration at node T is first computed through coordinate transformation. Subsequently, reaction forces and moments acting at the two attachment points are calculated using Eqs. 36-38. These loads are passed to the structural dynamics modules through the glue code interface, enabling fully coupled aero-hydro-servo-elastic analysis of TMDI-equipped wind turbines.

The above modifications are summarized in terms of a flow chart in Fig. A-1 and these integrate the TMDI configuration within the existing StC environment with minimal intrusion to the established module architecture. The implementation maintains full backward compatibility with existing TMD configurations and other structural control modes supported by OpenFAST. Furthermore, the dynamic node allocation strategy ensures computational efficiency by incurring additional overhead only when TMDI functionality is explicitly activated.

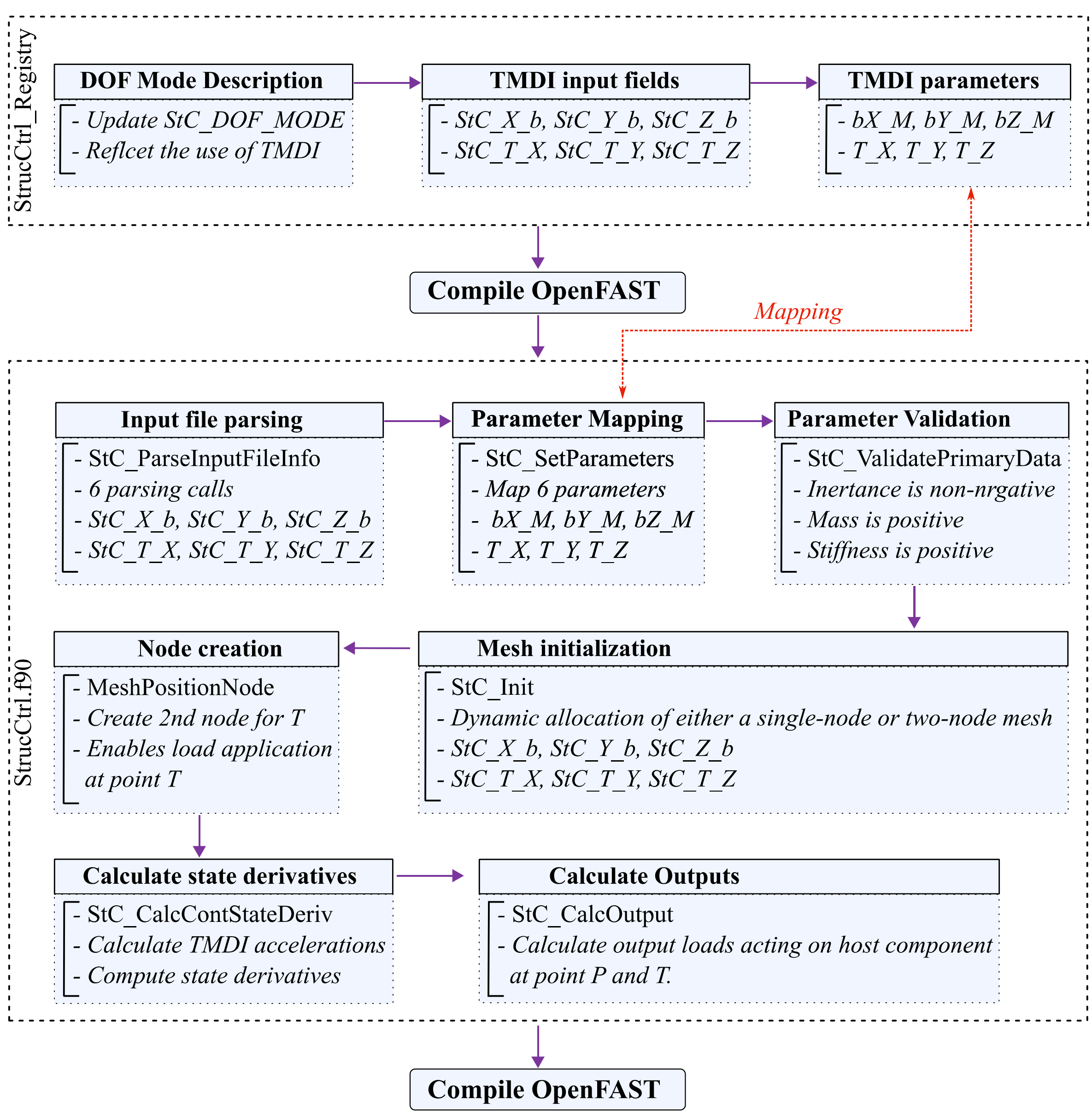

**Fig. A-1** Flowchart illustrating the TMDI integration workflow within OpenFAST Structural Control (StC) module, showing registry file modifications and computational subroutine modifications.

## Appendix B- Design load cases for IEA 15MW reference wind turbine

The metocean conditions considered to specify design load cases (DLCs) in Sections 3 and 4 for the IEA 15MW reference wind turbine as per as the IEC 61400-1 stem from wind and wave measurements acquired over several years in the SEM-REV offshore testing site [47], located 20 km off the coast of Le Croisic, France (47°14.34' N, 2°46.08' W).

For the extreme wind event with 50 years of return period the extrapolated SEM-REV mean wind speed at hub height, $\bar{U}_{\text{hub}}$, (150 m above MSL) is 40.55 m/s. Further, the turbulent wind components at hub height in the fore-aft, side-side, and vertical directions (x, y, z) are specified in OpenFAST using the following one-sided Kaimal-type spectra, prescribed by IEC 61400-1 [60]

$$S_x(f) = \frac{4\sigma_u^2\left(\frac{340.2}{\bar{U}_{hub}}\right)}{\left[1+6f\left(\frac{340.2}{\bar{U}_{hub}}\right)\right]^{5/3}}, \qquad S_y(f) = \frac{4(0.8\sigma_u)^2\left(\frac{113.4}{\bar{U}_{hub}}\right)}{\left[1+6f\left(\frac{113.4}{\bar{U}_{hub}}\right)\right]^{5/3}}, \qquad S_z(f) = \frac{4(0.5\sigma_u)^2\left(\frac{27.72}{\bar{U}_{hub}}\right)}{\left[1+6f\left(\frac{27.72}{\bar{U}_{hub}}\right)\right]^{5/3}}, \tag{B.1}$$

respectively. In the above expressions, $f$ is the frequency in Hz and $\sigma_u$ is the standard deviation of the turbulent wind component in the fore-aft direction at hub height, given as $\sigma_u = 0.096(0.75\bar{U}_{\text{hub}} + 5.6\ m/s)$ for the normal turbulence model with low turbulence characteristics, or as $\sigma_u = 0.11\bar{U}_{\text{hub}}$ for the extreme turbulence model per IEC 61400-1 [60].

For sea wave modelling the one-sided JONSWAP-like spectrum specified by IEC 61400–3 [46] as

$$S_w(f) = \frac{5}{16}\frac{H_s^2 f_p^4}{f^5}\exp\left[-\frac{5}{4}\left(\frac{f_p}{f}\right)^4\right](1-0.287\ln\gamma)\gamma^a, \tag{B.2}$$

where $H_s$ is the significant wave height, $f_p$ is the peak wave frequency, and $\gamma$ is the peak enhancement factor given as defined as

$$\gamma = \begin{cases} 5, & \text{for } \dfrac{T_p}{\sqrt{H_s}} \leq 3.6, \\ \exp\left(5.75 - 1.15\dfrac{T_p}{\sqrt{H_s}}\right), & \text{for } 3.6 \leq \dfrac{T_p}{\sqrt{H_s}} \leq 5, \\ 1, & \text{for } \dfrac{T_p}{\sqrt{H_s}} > 5, \end{cases} \tag{B.3}$$

where $T_p = 1/f_p$ is the peak wave period. Further, in Eq. (B.2), the exponent $a$ is given as

$$a = \exp\left[-\frac{(f - f_p)^2}{2\sigma^2 f_p^2}\right], \tag{B.4}$$

where the shape parameter $\sigma$ is 0.07 for $f \leq f_p$ or 0.09 for $f > f_p$. For the SEM-REV site, significant wave height of 9.25 m for 50-year return period is found with corresponding peak enhancement factor of 1.708 based on wave height field measurements taken almost continuously in the period of 2009-2022. Moreover, SEM-REV site-specific operational wave conditions are represented by $H_s = 0.75$ m, $T_p = 10.5$ s, and $\gamma = 1$ based on the annual directional sample distribution for the period of 1994 - 2019.

Overall, Table B.1 summarizes metocean conditions and IEC 61400-specific modelling assumptions which have been used in Sections 3 and 4 of this paper to define wind and wave DLCs in OpenFAST.

**Table B.1** Specification of design load cases (DLCs) per IEC 61400 for the SEM-REV testing site [47].

| **Load case designation** | **1.2** | **1.6** | **6.1** |
|---|---|---|---|
| **State of turbine** | Power production | Power production | Parked (idling) |
| **Limit state** | Fatigue | Ultimate | Ultimate |
| **Wind turbulence model** | Normal | Normal | Extreme |
| **Wave conditions** | Operational | Extreme | Extreme |
| **Return period** | very short* | 50 years | 50 years |
| **Mean wind speed at hub height, $\overline{U}_{hub}$ (m/s)** | 10.59** | 10.59** | 40.55 |
| **Wave height $H_s$ (m)** | 0.75 | 9.25 | 9.25 |
| **Peak wave period $T_p$ (s)** | 10.5 | 15.79 | 15.79 |
| **Gamma ($\gamma$)** | 1.00 | 1.708 | 1.708 |

* Denotes frequently encountered environmental loads (operational state)

** Rated wind speed of the IEA 15MW turbine [44]